\documentclass[mathematics,article,accept,pdftex,oneauthor]{Definitions/mdpi} 

\firstpage{1} 
\pubvolume{1}
\issuenum{1}
\articlenumber{0}
\pubyear{2024}
\copyrightyear{2024}
\externaleditor{Academic Editor: Francesco Aldo Costabile }
\datereceived{31 July 2024} 
\daterevised{31 August 2024} 
\dateaccepted{4 September 2024} 
\datepublished{ } 
\hreflink{https://doi.org/} 

\usepackage{mathrsfs}
\usepackage{tikz}
\usepackage{amsmath, amsthm, amsbsy, amssymb}

\Title{On Polar Jacobi Polynomials}

\TitleCitation{On Polar Jacobi Polynomials}

\Author{{Roberto} 
 S. Costas-Santos \orcidA{}}

\AuthorNames{Roberto S. Costas-Santos}

\AuthorCitation{{Costas-Santos}
, R. S.}

\address[1]{Department of Quantitative Methods, Universidad Loyola Andaluc\'ia, E-41704 Seville, Spain; rscosa@gmail.com
}

\abstract{In the present work, we investigate certain algebraic and differential
properties of the orthogonal polynomials with respect to a 
discrete--continuous Sobolev-type inner product defined in terms of 
the Jacobi measure.}


\newcommand\rcrfi[1]{{\color{black}{#1}}} 

\begin{document}
\section{Introduction}
Let us consider a measure $\mu$ supported on the subset of the 
complex plane $\gamma$. 

\rcrfi{In the vector space of polynomials with complex coefficients $\mathbb P$, 
let us  introduce the following inner product:}
\begin{equation} \label{eq:orthsob}
\langle f,g\rangle_\xi:=M f(\xi)g(\xi)
+\int_\gamma f'(z) g'(z) d\mu(z),\quad M, \xi\in \mathbb C,
\end{equation}
assuming that the integral exists. The inner product 
{equation} 
 \eqref{eq:orthsob} is called 
\rcrfi{} 
discrete--continuous Sobolev-type, 
which is a particular case of the Sobolev-type inner products. 
Algebraic and analytical properties and the \rcrfi{limiting} behavior 
of the families of orthogonal polynomials with respect to 
Sobolev-type inner products \rcrfi{have been thoroughly} used
for the last 25 years.
For an overview of this subject, {see} 
 \cite{mr1808577}, or the 
introduction of \cite{MR2271726} as well as 
\cite{mr3360352}.

\rcrfi{Discrete--continuous Sobolev inner products were} introduced 
in \cite{MR294963} to study the behavior of the \rcrfi{optimal} polynomial 
approximation of absolutely continuous functions in the norm 
generated by a Sobolev inner product as {equation} \eqref{eq:orthsob}. 
Later, in \cite{MR1080667}, R. Koekoek considered the 
Laguerre case with $d\mu = x^\alpha e^{-x} dx$, $\alpha>-1$, 
$\gamma = [0,+\infty)$ and 
$\xi = 0$. The Gegenbauer case was studied by Bavinck and 
Meijer in \cite{MR1013456,MR1054761}, with 
$d\mu=(1-x^2)^{\lambda-1/2}dx$, 
$\lambda>-1/2$, $\gamma = [-1,1]$ and $\xi_1 =-1$ and 
$\xi_2 = 1$.

The families of orthogonal polynomials with respect to 
{equation} \rcrfi{\eqref{eq:orthsob}} have been studied as \rcrfi{an} 
extension of the Bochner--Krall theory (i.e., 
\rcrfi{polynomial sequences which} are simultaneously eigenfunctions of a differential 
operator and orthogonal with respect to an inner product\rcrfi{), }
see for the discrete--continuous 
\rcrfi{case}
\cite{mr1795525, MR1641942,mr1651585}\rcrfi{.}

The starting point of our work is the orthogonality with 
respect to the Jacobi case. 
Let $(Q_n)$ be the 
\rcrfi{monic  orthogonal polynomial sequence} with respect to the 
Sobolev inner product
\begin{equation} \label{eq:orthsobjac}
\langle f,g\rangle:= f(\xi)g(\xi)
+\int_\gamma f'(z) g'(z) (1-z)^\alpha (1+z)^\beta dz,
\end{equation}
where $\alpha, \beta \in \mathbb C$, $\gamma$ is a 
path encircling the points $+1$ and $-1$ first in a positive 
sense and then in a negative 
sense, as shown in {Figure 2.1} 
 in \cite{mr2149265}, $M = 1$ and 
$\xi \in \mathbb C$.

Let us consider $\left(P_n^{(\alpha,\beta)}(z)\right)$ as the 
monic Jacobi polynomials, i.e., the 
\rcrfi{monic orthogonal polynomial sequence} 
with  respect to the measure \rcrfi{$d\mu(z)=(1-z)^\alpha (1+z)^\beta dz$}. 
Therefore, for each $n$\rcrfi{$\in \mathbb N$}, 
$P_n^{(\alpha,\beta)}(z)$ satisfies
\begin{equation} \label{eq:orthjac}
\int_\gamma P_n^{(\alpha,\beta)}(z)z^k (1-z)^\alpha 
(1+z)^\beta dz=0,
\quad k= 0, 1, \ldots, n-1,
\end{equation}
and let $\Pi_{n+1,\xi}$ be the polynomial primitive of 
$(n+1)P_n^{(\alpha,\beta)}$ \rcrfi{that has a zero} at $\xi$, i.e., 
for all $n\ge 1$, we have
\begin{equation} \label{eq:primitive}
\Pi_{n+1,\xi}(\xi)=0,\qquad \frac{d}{dz}\Pi_{n+1,\xi}(z)=
(n+1)P_n^{(\alpha,\beta)}(z).
\end{equation}
{Then,} 
 by definition of $Q_n$, it is clear that 
$\Pi_{n+1,\xi}(z) = Q_{n+1}(z)$ 
for all $n=0,1, 2, \ldots$. It is straightforward to prove 
that the 
\rcrfi{polynomial of degree $n$ which is orthogonal} 
with respect to {equation} \eqref{eq:orthsobjac} can be written 
as follows ({{one must} 
 consider $Q_0(z)=1$}): 
\[
Q_n(z)=(z-\xi)P_{n-1}(z),
\]
where $P_n$ is called {the} 
{\tt polar polynomial associated 
with} $\mu$ 
(see \cite{MR2764236}) and $\xi$ from now on will be called 
the {\tt pole}. 
Let us define the differential operator 
$L_\xi : {\cal H}^1(\gamma)\rightarrow 
L^2(\gamma)$ as
\begin{equation} \label{eq:operatorL}
L_\xi[f(z)]=f(z)+(z-\xi) \frac{d}{dz}f(z),
\end{equation} 
where $ {\cal H}^1(\gamma):=\{f\in 
L^2(\gamma): f'(z)\in L^2(\gamma)\}$ 
is the Sobolev space of index 1. Taking into account that
$Q_n$ is orthogonal with respect to the inner 
product {equation} \eqref{eq:orthsobjac}, we have
\[
\int_\gamma L_\xi[P_n(z)]z^k (1-z)^\alpha (1+z)^\beta dz=
\int_\gamma \left(P_n(z)+(z-\xi)P'_n(z)\right)z^k (1-z)^\alpha 
(1+z)^\beta dz=0,
\]
for $k = 0, 1, \ldots, n-1$. Therefore, $P_n$ is the 
\rcrfi{monic orthogonal polynomial of degree $n$} 
with respect to the differential operator 
$L_\xi$, and the measure
\rcrfi{$d\mu$ 
(see \cite{MR1934900, MR3045189, MR3351532, MR3032627, MR2764236, MR4089668})}. 
In such a case, 
\rcrfi{for $n\in \mathbb N_0$,}
we have
\begin{equation} \label{eq:relpoljac}
P_n(z)+(z-\xi)P'_n(z)=(n+1)P_n^{(\alpha,\beta)}(z).
\end{equation} 

\rcrfi{The main aim of this work is to study algebraic (zero 
localization) and  differential  properties of 
the polynomial sequences that are orthogonal with respect 
to the inner product {equation} \eqref{eq:orthsobjac} for the Jacobi weight case, 
which is a natural  extension of the Legendre case \cite{MR2764236}}.

In Section \ref{sec:2}, 
we obtain several algebraic relations between the polar Jacobi 
polynomials and the Jacobi polynomials and some differential 
and different identities related to the polar Jacobi polynomials. 
\rcrfi{Finally, in} Section \ref{sec:3} we study the location of the zeros for the 
polynomials $P_n$. 
\section{Algebraic Properties of the Polar Jacobi Polynomials}
\label{sec:2}
Let us start by summarizing some basic properties of the Jacobi 
orthogonal polynomials to be used in the sequel to {Chapter 18} in \cite{DLMF}

\begin{Proposition}
Let $\left(P_n^{(\alpha,\beta)}(z)\right)$ be the classical 
monic Jacobi orthogonal polynomial sequence.
The following statements hold:
\begin{enumerate}
\item Three-term recurrence relation.
\begin{equation} \label{eq:JTTRR}
P_{n+1}^{(\alpha,\beta)}(z)=(z-\beta_n) 
P_{n}^{(\alpha,\beta)}(z)-
\gamma_n P_{n-1}^{(\alpha,\beta)}(z),\quad n=0, 1, \ldots,
\end{equation}
with initial condition $P_{0}^{(\alpha,\beta)}(z)=1$, and 
recurrence coefficients

\begin{adjustwidth}{-\extralength}{0cm}
\centering 
\[
\beta_n=\dfrac{\beta^2-\alpha^2}
{(\alpha+\beta+2n)(\alpha+\beta+2n+2)},
\quad
\gamma_n=\dfrac{4n(\alpha+n)(\beta+n)(\alpha+\beta+n)}
{(\alpha+\beta+2n-1)(\alpha+\beta+2n)^2(\alpha+\beta+2n+1)}.
\]
\end{adjustwidth}
\item First structure relation.
\begin{equation} \label{eq:JFSR}
(1-z^2)\frac{d}{dz} P_{n}^{(\alpha,\beta)}(z)=
-n P_{n+1}^{(\alpha,\beta)}(z)+
\widehat\beta_n P_{n}^{(\alpha,\beta)}(z)
+\widehat\gamma_nP_{n-1}^{(\alpha,\beta)}(z),\quad n=0, 1, \ldots,
\end{equation}
with coefficients 
\vspace{-6pt}{}
\begin{adjustwidth}{-\extralength}{0cm}
\centering 
\[
\widehat\beta_n=\dfrac{2n(\alpha-\beta)(\alpha+\beta+n+1)}
{(\alpha+\beta+2n)(\alpha+\beta+2n+2)},
\quad
\widehat\gamma_n=\dfrac{4n(\alpha+n)(\beta+n)
(\alpha+\beta+n)(\alpha+\beta+n+1)}
{(\alpha+\beta+2n-1)(\alpha+\beta+2n)^2(\alpha+\beta+2n+1)}.
\]
\end{adjustwidth}
\item Second structure relation.
\begin{equation} \label{eq:JSSR}
P_{n+1}^{(\alpha-1,\beta-1)}(z)=P_{n+1}^{(\alpha,\beta)}(z)+
\widetilde\beta_n P_{n}^{(\alpha,\beta)}(z)
+\widetilde\gamma_nP_{n-1}^{(\alpha,\beta)}(z),\quad n=0, 1, \ldots,
\end{equation}
\vspace{-24pt}{}
\begin{adjustwidth}{-\extralength}{0cm}
\centering 
\[
\widetilde\beta_n=\dfrac{(2n+2)(\alpha-\beta)}
{(\alpha+\beta+2n)(\alpha+\beta+2n+2)},
\quad
\widetilde\gamma_n=-\dfrac{4n(n+1)(\alpha+n)(\beta+n)}
{(\alpha+\beta+2n-1)(\alpha+\beta+2n)^2(\alpha+\beta+2n+1)}.
\]
\end{adjustwidth}
\item Squared Norm. For every $n\ge 0$,
\vspace{-12pt}{}
\begin{adjustwidth}{-\extralength}{0cm}
\centering 
\begin{equation} \label{eq:JSN}
\left\|P_{n}^{(\alpha,\beta)}(z)\right\|^2 = \int_\gamma 
\left(P_{n}^{(\alpha,\beta)}(z)\right)^2 
\omega(z;\alpha,\beta)\, dz=\dfrac{\rcrfi{2^{2n+\alpha+\beta+1}}n! 
\Gamma(\alpha+n+1)\Gamma(\beta+n+1)
\Gamma(\alpha+\beta+n+1)}
{\Gamma(\alpha+\beta+2n+1)\Gamma(\alpha+\beta+2n+2)}.
\end{equation}
\end{adjustwidth}
\item Second-order difference equation. For every $n\ge 0$, 
\vspace{-12pt}{}
\begin{adjustwidth}{-\extralength}{0cm}
\centering 
\begin{equation} \label{eq:JSODE}
(1-z^2)\frac{d^2}{dz^2} P_{n}^{(\alpha,\beta)}(z)+
\left(\beta-\alpha-z(\rcrfi{\alpha+\beta+2})\right)
\frac{d}{dz}P_{n}^{(\alpha,\beta)}(z)=-n (\alpha+\beta+n+1) 
P_{n}^{(\alpha,\beta)}(z).
\end{equation}
\end{adjustwidth}
\item Forward shift operator.
\begin{equation} \label{eq:JFSO}
\frac{d}{dz}P_{n}^{(\alpha,\beta)}(z)= n 
P_{n-1}^{(\alpha+1,\beta+1)}(z), \quad n=0, 1, \ldots,
\end{equation}
\item Asymptotic formula. Let $z\in\mathbb C\setminus[-1,1]$. 
Put $\varphi(z)=z+\sqrt{z^2-1}$ where the branch of the square 
root is chosen so that $|z+\sqrt{z^2-1}|>1$ for $z\in 
\mathbb C\setminus [-1,1]$. Then,
\begin{equation} \label{eq:JAF}
P_{n}^{(\alpha,\beta)}(z)=\dfrac{\varphi^n(z)}
{\sqrt n}\left(c(\alpha,\beta,z)
+{\cal O}(n^{-1})\right),
\end{equation} 
where $c(\alpha,\beta, z)$ is a function of $\alpha$ and $\beta$ 
and $x$ independent of $n$. The relation holds uniformly 
on compact sets of $\mathbb C\setminus[-1,1]$.
\end{enumerate}
\end{Proposition} 

Let us obtain the algebraic relations between the Jacobi 
polynomials and the polar Jacobi polynomials.
\begin{Lemma}\label{lem:1} For any $\alpha, \beta, \xi\in
\mathbb C$. 
The polar Jacobi polynomials can be written in terms of the 
Jacobi polynomials as follows:

\begin{adjustwidth}{-\extralength}{0cm}
\centering 
\begin{eqnarray}
\label{eq:RJPoJ2}P_n(z)&=&\dfrac
{P_{n+1}^{(\alpha-1,\beta-1)}(z)-
P_{n+1}^{(\alpha-1,\beta-1)}(\xi)}{z-\xi}\\
\nonumber &=&\frac 1{\alpha+\beta+n}\left[(z+\xi)
\frac{d}{dz}P_{n}^{(\alpha,\beta)}(z)+
\dfrac{\xi^2-1}{z-\xi}\left(\frac{d}{dz}P_{n}^{(\alpha,\beta)}(z)
-\frac{d}{dz}P_{n}^{(\alpha,\beta)}(\xi)\right)\right.\\
\label{eq:RJPoJ3}&&\left.+(\alpha+\beta)
P_{n}^{(\alpha,\beta)}(z)
+\dfrac{\alpha-\beta+\xi(\alpha+\beta)}{z-\xi}
\left(P_{n}^{(\alpha,\beta)}(z)-
P_{n}^{(\alpha,\beta)}(\xi)\right)\right].
\end{eqnarray}
\end{adjustwidth}
{Therefore,} 
\begin{equation} \label{eq:RJPoJ1}
(z-\xi)P_n(z)=P_{n+1}^{(\alpha,\beta)}(z)+
\widetilde\beta_{n} P_{n}^{(\alpha,\beta)}(z)
+\widetilde\gamma_{n}P_{n-1}^{(\alpha,\beta)}(z)-
P_{n+1}^{(\alpha-1,\beta-1)}(\xi).
\end{equation}
\end{Lemma}
\begin{proof} 
From {equation} \eqref{eq:relpoljac}, we have 
\[
(n+1)P_{n}^{(\alpha,\beta)}(z)= \dfrac{d}{dz}\left((z-\xi)
P_n(z)\right). 
\]
Therefore, by using the forward shift operator {equation} \eqref{eq:JFSO}, 
we have
\begin{equation} \label{eq:relJPoint}
(z-\xi)P_n(z)=(n+1)\int_\xi^z P_{n}^{(\alpha,\beta)}(z)\, dz= 
\int_\xi^z \dfrac{d}{dz} 
\left(P_{n+1}^{(\alpha-1,\beta-1)}(z)\right)\, dz.
\end{equation} 
{From} {equation} \eqref{eq:relJPoint}, the identity 
\eqref{eq:RJPoJ2} follows.
\rcrfi{Using} 
{equation} \eqref{eq:RJPoJ2} 
and 
the second structure relation 
\eqref{eq:JSSR} the 
expression 
\eqref{eq:RJPoJ1} follows. 

Let us now to prove {equation} \eqref{eq:RJPoJ3}:
by using the second-order differential 
\rcrfi{Equation \eqref{eq:JSODE}} 
\rcrfi{with $n, \alpha,\beta$ being shifted to $n+1,\alpha-1,\beta-1$, 
respectively,}
and the forward 
shift operator of the Jacobi polynomials, we obtain
\begin{equation} \label{eq:aux1}
(1-z^2)\frac{d}{dz} P_{n}^{(\alpha,\beta)}(z)+
\left(\beta-\alpha-z(\alpha+\beta)\right)
P_{n}^{(\alpha,\beta)}(z)=-(\alpha+\beta+n) 
P_{n+1}^{(\alpha-1,\beta-1)}(z).
\end{equation}
{By} using the differential Equation \eqref{eq:aux1} and the identity
\[
\dfrac{f(z)g(z)-f(\xi)g(\xi)}{z-\xi}=
\dfrac{f(z)-f(\xi)}{z-\xi}g(z)+f(\xi)\dfrac{g(z)-g(\xi)}{z-\xi},
\]
then {equation} \eqref{eq:RJPoJ3} follows and hence the result holds.
\end{proof}
The following additional property of orthogonality holds.
\begin{Theorem} \label{thm:1}
The polar Jacobi polynomial $P_n$ with pole $\xi\in \mathbb C$ 
fulfills the following property of orthogonality:
\vspace{-12pt}{}
\begin{adjustwidth}{-\extralength}{0cm}
\centering 
\begin{equation} \label{eq:JPoOrt}
\int_\gamma \left(P_n(z)+(z-\xi)\frac{d}{dz}P_n(z)\right)
P_{m}^{(\alpha,\beta)}(z)
\omega(z;\alpha,\beta) dz=
\left\{\begin{array}{ll}0,&m\ne n,\\[3mm]
(n+1)\left\|P_{n}^{(\alpha,\beta)}(z)\right\|^2,&m=n.
\end{array}\right.
\end{equation}
\end{adjustwidth}
{Furthermore}, if $n>1$, then 

\begin{adjustwidth}{-\extralength}{0cm}
\centering 
\begin{equation} \label{eq:JPoint2}
\int_\gamma(z-\xi)P_n(z)P_{m}^{(\alpha,\beta)}(z)
\omega(z;\alpha,\beta)dz=\left\{
\begin{array}{ll}
-\dfrac{\rcrfi{2^{\alpha+\beta+1}}\Gamma(\alpha+1)\Gamma(\beta+1)}
{\Gamma(\alpha+\beta+2)}
P_{n+1}^{(\alpha-1,\beta-1)}(\xi),&m=0,\\[2mm]
0,&0<m<n-1,\\[2mm]
\widetilde\gamma_n \left\|P_{n-1}^{(\alpha,\beta)}(z)
\right\|^2, &n-1=m,\\[2mm]
\widetilde\beta_n \left\|P_{n}^{(\alpha,\beta)}(z)\right\|^2, 
&n=m,\\[2mm]
\left\|P_{n+1}^{(\alpha,\beta)}(z)\right\|^2, &n+1=m,\\[2mm]
0,&n+1<m.
\end{array}\right.
\end{equation} 
\end{adjustwidth}
\end{Theorem}
\begin{proof} 
Taking into account {equation} \eqref{eq:relpoljac}, we have 
\vspace{-12pt}{}
\begin{adjustwidth}{-\extralength}{0cm}
\centering 
\[
\int_\gamma \left(P_n(z)+(z-\xi)\frac d{dz}P_n(z)
\right)P_m^{(\alpha,\beta)}(z) 
\omega(z;\alpha,\beta)\, dz=(n+1)\int_\gamma 
P_n^{(\alpha,\beta)}(z) 
P_m^{(\alpha,\beta)}(z) \omega(z;\alpha,\beta)\, dz.
\]
\end{adjustwidth}
{So}, the first property of orthogonality follows.
By using the relation 
\eqref{eq:RJPoJ1} and considering 
the property of orthogonality of the Jacobi polynomials, 
the second property of orthogonality follows. Hence, the 
result holds.
\end{proof}

\begin{Theorem} \label{thm:2}
The sequence of polar Jacobi polynomials $(P_n)$ with pole 
$\xi\in \mathbb C$ satisfies the following recurrence relation:
\begin{equation} \label{eq:PoJRR} 
P_{n+1}(z) = zP_n(z) + a_nP_n(z) + b_nP_{n-1}(z)+
P_{n+1}^{(\alpha-1,\beta-1)}(\xi),\quad 
n=0, 1, \ldots, 
\end{equation} 
with initial conditions $P_{-1}(z) = 0$ and $P_0(z) = 1$, 
and coefficients
\vspace{-12pt}{}
\begin{adjustwidth}{-\extralength}{0cm}
\centering 
\[
a_n=\dfrac{(\alpha+\beta-2)(\alpha-\beta)}{(\alpha+\beta+2n)(\alpha+\beta+2n+2)},
\quad 
b_n=-\dfrac{4(n+1)(\alpha+n)(\beta+n)(\alpha+\beta+n-1)}
{(\alpha+\beta+2n-1)(\alpha+\beta+2n)^2(\alpha+\beta+2n+1)}.
\]
\end{adjustwidth}
\end{Theorem}
\begin{proof} 
Let the sequence $(\nu_{n,k})$ be such that
\[
(z-\xi)P_n(z)=\sum_{k=0}^{n+1} \nu_{n,k} P_k(z),
\]
\rcrfi{where $\nu_{n,n+1}=1$.}

Then, by using {equation} \eqref{eq:relpoljac}, we obtain 
\begin{eqnarray} \nonumber
(z-\xi)\left(P_n(z)+(n+1)P_n^{(\alpha,\beta)}(z)\right)&=&
(z-\xi)\left(P_n(z)+((z-\xi)P_n(z))'\right)\\[3mm]
&=&\label{eq:thm2-1}\sum_{k=0}^{n+1} \nu_{n,k} \left(P_k(z)+(z-\xi)
\frac d{dz}P_k(z)\right).
\end{eqnarray}
{By} the property of orthogonality {equation} \eqref{eq:JPoOrt}, we have
\vspace{-12pt}{}
\begin{adjustwidth}{-\extralength}{0cm}
\centering 
\begin{equation} \label{eq:thm2-2}
\sum_{k=0}^{n+1} 
\rcrfi{\nu_{n,k}}
\int_\gamma \left(P_k(z)+(z-\xi)
\frac d{dz}P_k(z)\right)P_m^{(\alpha,\beta)}(z) 
\omega(z;\alpha,\beta)\, dz=
\rcrfi{\nu_{n,k}}
 (m+1) 
\left\|P_m^{(\alpha,\beta)}(z)\right\|^2,
\end{equation}
\end{adjustwidth}
for $m= 0, 1, \ldots, n$.

On the other hand, let us denote
\[
I_{n,m}=\int_\gamma (z-\xi)P_m^{(\alpha,\beta)}(z) \left(P_n(z)+(n+1)
P_n^{(\alpha,\beta)}(z)\right)
\omega(z;\alpha,\beta)\, dz.
\]
{From} the orthogonality of the Jacobi polynomials and the property of 
orthogonality 
\eqref{eq:JPoint2}, we obtain
\begin{equation} \label{eq:thm2-3}
I_{n,m}=\left\{\begin{array}{ll}
-\rcrfi{2^{\alpha+\beta+1}}P_{n+1}^{(\alpha-1,\beta-1)}(\xi)\dfrac{\Gamma(\alpha+1)\Gamma(\beta+1)}
{\Gamma(\alpha+\beta+2)},&m=0,\\[2mm]
0,& 0<m<n-1,\\[2mm]
-\widetilde \gamma_n (\alpha+\beta+n-1)\left\|P_{n-1}^{(\alpha,\beta)}(z)
\right\|^2,&n-1=m\\[2mm]
\left(\dfrac{2-\alpha-\beta}{2}\widetilde \beta_n-\xi(n+1)\right)
\left\|P_{n}^{(\alpha,\beta)}(z)
\right\|^2,& n=m.
\end{array} \right.
\end{equation}
{Thus}, multiplying {equation} \eqref{eq:thm2-1} by $P_m^{(\alpha,\beta)}(z)$, 
integrating over $\gamma$, and using {equations} \eqref{eq:thm2-2} and 
\eqref{eq:thm2-3}, we obtain
\[
\nu_{n,m}=\dfrac{I_{n,m}}{(m+1)
\left\|P_{m}^{(\alpha,\beta)}(z)\right\|^2}
=\left\{\begin{array}{ll}
-P_{n+1}^{(\alpha-1,\beta-1)}(\xi),&m=0,\\[2mm]
0,&0<m<n-1,\\[2mm]
-b_n,&m=n-1,\\[2mm]
\dfrac{2-\alpha-\beta}{2(n+1)}\widetilde \beta_n-\xi, &m=n.
\end{array} \right.
\]
{The} expression 
\eqref{eq:PoJRR} is obtained after a 
straightforward calculation.
\end{proof}

A direct consequence of this result is the following.
\begin{Corollary}[The polar ultraspherical case] \label{cor:1}
The sequence of symmetric polar Jacobi polynomials with pole 
$\xi\in \mathbb C$, i.e., the sequence of polar ultraspherical polynomial with 
pole $\xi\in \mathbb C$, satisfies, namely {${\sf P}_n$}
, the following recurrence relation:
\vspace{-12pt}{}
\begin{adjustwidth}{-\extralength}{0cm}
\centering 
\begin{equation} \label{eq:PogePRR} 
{\sf P}_{n+1}(z) = z{\sf P}_n(z) - \dfrac{(n+1)(2\alpha+n-1)}
{(2\alpha+2n-1)(2\alpha+2n+1)}{\sf P}_{n-1}(z)+
{\sf P}_{n+1}^{(\alpha-1,\alpha-1)}(\xi),\quad 
 n=0, 1, \ldots, 
\end{equation}
\end{adjustwidth} 
with initial conditions ${\sf P}_{-1}(z) = 0$ and ${\sf P}_0(z) = 1$.
\end{Corollary}

Another direct consequence is the fact 
\rcrfi{that} 
when one, or both, of the parameters is a negative integer, 
we can factorize 
the Jacobi polynomial. In fact \rcrfi{(see {(1.2) in} \cite{mr2149265})}, 
\begin{align}
P_n^{(-k,\beta)}(z)=&(z-1)^k P^{(k,\beta)}_{n-k},\label{eq:26}\\
P_n^{(\alpha,-k)}(z)=&(z+1)^k P^{(\alpha,k)}_{n-k}\label{eq:27}.
\end{align}

\begin{Remark}\label{rem:1}
Since in some results we will consider the polar Jacobi polynomials
with different parameters and poles, to avoid such possible confusion, 
we will denote by {$P_n(z;\alpha,\beta;\xi)$} the polar Jacobi polynomial 
of degree $n$ with parameters $\alpha$ and $\beta$, and pole at $\xi$.
\end{Remark}

\begin{Corollary}[The factorization] \label{cor:2}
For any positive integer $k$, the following identities hold:
\begin{align}
P_{n+k}(z;-k,\beta;1)=&(z-1)^k P^{k+1,\beta-1}_{n}(z)\label{eq:28}\\
=&(z-1)^k \left((z-1)P_{n-1}(z;k+2,\beta;1)
+ P^{k+1,\beta-1}_{n}(1)\right),\label{eq:29}\\
P_{n+k}(z;\alpha,-k;-1)=&(z+1)^k P^{\alpha-1,k+1}_{n}(z)\label{eq:30}\\
=&(z+1)^k \left((z+1)P_{n-1}(z;\alpha,k+2;-1)
+ P^{\alpha-1, k+1}_{n}(-1)\right).\label{eq:31}
\end{align}
{Moreover}, the recurrence coefficients satisfy the following relations:
\[
a_{n+k}(-k,\beta;1)=a_{n-1}(k+2,\beta;1),\quad 
b_{n+k}(-k,\beta;1)=b_{n-1}(k+2,\beta;1),
\]
and
\[
a_{n+k}(\alpha,-k;-1)=a_{n-1}(\alpha,k+2;-1),\quad 
b_{n+k}(\alpha,-k;-1)=b_{n-1}(\alpha,k+2;-1).
\]
\end{Corollary}

\begin{proof}
The identities 
\eqref{eq:28}--\eqref{eq:31} 
follow by using the factorization of the Jacobi polynomials {equations} \eqref{eq:26} and 
\eqref{eq:27}. In order to obtain the relation between the recurrence
 coefficients defined in Theorem \ref{thm:2}, we must use the former 
 factorization(s), 
 and after a straightforward calculation the identities 
 follow.
\end{proof} 
The last result of this section is due the parity relation of the 
Jacobi polynomials, i.e.,
\begin{equation} \label{eq:33}
P_n^{(\alpha,\beta)}(z)=
\rcrfi{(-1)^n} 
P_n^{(\beta,\alpha)}(-z).
\end{equation}
\begin{Lemma} \label{lem: 2}
For any $\xi\in \mathbb C$, the following identity holds:
\begin{equation} \label{eq:34}
P_n(z;\alpha,\beta;\xi)=(-1)^nP_n(-z;\beta,\alpha;-\xi).
\end{equation}
\end{Lemma}

\begin{proof} 
Starting from 
\rcrfi{{equation} \eqref{eq:RJPoJ2}} 
and using {equation} \eqref{eq:33}, we have
\begin{align*}
P_n(z;\alpha,\beta;\xi)=&\dfrac{P^{\alpha-1,\beta-1}_{n+1}(z)
-P^{\alpha-1,\beta-1}_{n+1}(\xi)}{z-\xi}=
(-1)^{n+1}\dfrac{P^{\beta-1,\alpha-1}_{n+1}(-z)
-P^{\beta-1,\alpha-1}_{n+1}(-\xi)}{z-\xi}\\
=&(-1)^{n}\dfrac{P^{\beta-1,\alpha-1}_{n+1}(-z)
-P^{\beta-1,\alpha-1}_{n+1}(-\xi)}{-z-(-\xi)}=
(-1)^{n}P_n(-z;\beta,\alpha;-\xi).
\end{align*}
\end{proof}
\section{Zero Location}
\label{sec:3}
Finding the roots of polynomials is a problem of interest in both 
mathematics and in areas of application such as physical systems, 
which can be reduced to solving certain 
\rcrfi{equations}. 
There are very interesting geometric relationships between the 
roots of a polynomial $f_n(z)$ and those of $f'_n (z)$. 
The most important result is the following.

\begin{Theorem}[The Gau{\ss}--Lucas theorem \cite{zbMATH2717048}] \label{thm:3}
Let $f_n(z)\in \mathbb C[z]$ be a polynomial of degree of at least one. 
All zeros of $f_n'(z)$ lie in the convex hull of the zeros of the zeros 
of $f_n(z)$.
\end{Theorem}

In this section, we are going to study the zero distribution for the polar 
Jacobi polynomials. The next 
result, which was obtained by G. Szeg\H{o}, 
\rcrfi{is useful to }
estimate where such zeros are located.

\begin{Theorem}[Szeg\H{o}'s theorem \cite{MR1544526,MR1367960}]
\label{thm:4}
Let $a(z)$ and $b(z)$ be polynomials 
\rcrfi{of the form}
\[
a(z)=\sum_{\ell=0}^n a_\ell {n \choose \ell} z^\ell,\quad 
b(z)=\sum_{\ell=0}^n b_\ell {n \choose \ell} z^\ell. 
\]
\hl{If} the zeros of $a(z)$ lie in a closed disk $\overline D$ 
and $\lambda_1$, \ldots, $\lambda_n$ are the zeros of $b(z)$, then the 
zeros of the ``composition'' of the two
\[
c(z)=\sum_{\ell=0}^n a_\ell b_\ell {n \choose \ell} z^\ell, 
\]
have the form $\lambda_\ell \gamma_\ell$, where $\gamma_\ell\in 
\overline D$.
\end{Theorem}

By using this result, we are going to locate the 
\rcrfi{disk within which} 
all the zeros of the polar Jacobi are located. 

\begin{Theorem} \label{thm:5}
For any $\Re \alpha,\Re \beta>-1$ and $\xi\in \mathbb C$, the zeros of 
$P_n(z;\alpha,\beta;\xi)$ lie inside the closed disk 
$\overline{D}(0,2+|\xi|)$.
\end{Theorem}

\begin{proof} 
Starting from {equation} \eqref{eq:relpoljac} and assuming that
\[
a(z)=P_n^{(\alpha,\beta)}(w)=\sum_{k=0}^n \mu_k 
\rcrfi{w}^k,\quad
c(z)=P_n(z;\alpha,\beta;w
)=\sum_{k=0}^n \eta_k 
\rcrfi{w}^k,
\]
where $w:=z-\xi$, then $\eta_k=(n+1)/(k+1)\mu_k$. 
In order to apply Szeg\H{o}’s theorem, we consider
\[
b(z)=\sum_{k=0}^n \binom nk\dfrac{n+1}{k+1} w^k=
\sum_{k=0}^n \binom {n+1}{k+1}w^k=
\dfrac{(w+1)^{n+1}-1}{w}.
\]
{If} $b(w_1)=0$, then $|w_1+1|=1$, so $|z_1|\le 2+|\xi|$. Moreover, 
if $a(z_2)=0$\rcrfi{, i.e., $z_2$ is a zero of $a(z)$,} 
then $|z_2|\le 1$. Therefore, combining these inequalities 
and applying Szeg\H{o}'s theorem one obtains that if 
$c(z_3)=0$\rcrfi{, i.e., $z_3$ is a zero of $c(z)$,} 
 then $|z_3|\le 2+|\xi|$ and hence the result follows. 
\end{proof}

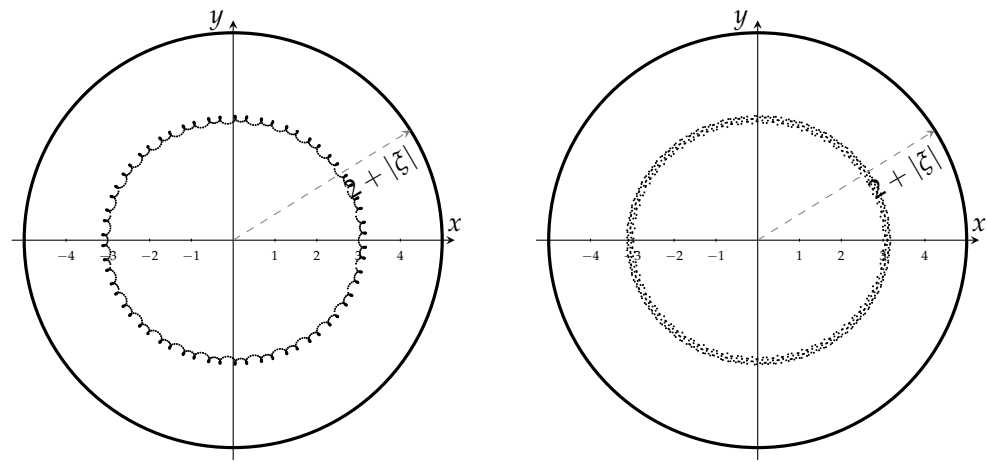
\begin{figure}[H]

\begin{minipage}{0.46\textwidth}
\begin{tikzpicture}[domain=-0.7:0.7,font=\sffamily, scale=0.55]
\draw[-stealth] (-5.3,0) -- (5.3,0) node[above] {$x$};
\draw[-stealth] (0,-5.3) -- (0,5.3) node[left] {$y$};
\draw plot [only marks, mark=*, mark options={fill=black},mark size=0.3pt]
coordinates{(-2.9376, -0.28396) (-2.9376, 
 0.28396) (-2.8163, -0.84028) (-2.8163, 
 0.84028) (-2.5787, -1.3622) (-2.5787, 
 1.3622) (-2.2344, -1.8285) (-2.2344, 
 1.8285) (-1.7977, -2.2199) (-1.7977, 
 2.2199) (-1.2865, -2.5206) (-1.2865, 
 2.5206) (-0.72164, -2.7180) (-0.72164, 
 2.7180) (-0.12631, -2.8042) (-0.12631, 
 2.8042) (0.47512, -2.7756) (0.47512, 
 2.7756) (1.0580, -2.6334) (1.0580, 
 2.6334) (1.5986, -2.3833) (1.5986, 
 2.3833) (2.0747, -2.0355) (2.0747, 
 2.0355) (2.4668, -1.6044) (2.4668, 
 1.6044) (2.7589, -1.1076) (2.7589, 
 1.1076) (2.9391, -0.56538) (2.9391, 0.56538) (-2.9516, 
 0.22729) (-2.9392, -0.34234) (-2.8421, 
 0.78762) (-2.8054, -0.89797) (-2.6153, 
 1.3157) (-2.5557, -1.4169) (-2.2803, 
 1.7901) (-2.2004, -1.8778) (-1.8510, 
 2.1911) (-1.7540, -2.2620) (-1.3450, 
 2.5026) (-1.2349, -2.5536) (-0.78295, 
 2.7116) (-0.66426, -2.7406) (-0.18792, 
 2.8096) (-0.065490, -2.8155) (0.41572, 
 2.7926) (0.53689, -2.7752) (1.0033, 
 2.6613) (1.1182, -2.6211) (1.5507, 
 2.4210) (1.6547, -2.3598) (2.0356, 
 2.0815) (2.1245, -2.0017) (2.4382, 
 1.6567) (2.5082, -1.5617) (2.7419, 
 1.1641) (2.7903, -1.0577) (2.9591, -0.51039) (3.0078, 
 0.057860) (-2.9787, 0.17879) (-2.9554, -0.39554) (-2.8791, 
 0.74581) (-2.8100, -0.95370) (-2.6605, 
 1.2823) (-2.5485, -1.4729) (-2.3321, 
 1.7664) (-2.1817, -1.9318) (-1.9071, 
 2.1782) (-1.7246, -2.3117) (-1.4031, 
 2.5009) (-1.1959, -2.5970) (-0.84076, 
 2.7213) (-0.61737, -2.7760) (-0.24303, 
 2.8302) (-0.012599, -2.8413) (0.36557, 
 2.8233) (0.59362, -2.7904) (0.96015, 
 2.7008) (1.1765, -2.6252) (1.5164, 
 2.4677) (1.7121, -2.3524) (2.0115, 
 2.1335) (2.1786, -1.9834) (2.4252, 
 1.7119) (2.5569, -1.5330) (2.8316, -1.0199) (2.9449, 
 0.67851) (2.9913, -0.46496) (3.0296, 0.10900) (-3.0140, 
 0.14284) (-2.9822, -0.43822) (-2.9214, 
 0.71806) (-2.8273, -1.0014) (-2.7081, 
 1.2639) (-2.5556, -1.5236) (-2.3829, 
 1.7581) (-2.1783, -1.9834) (-1.9592, 
 2.1803) (-1.7108, -2.3622) (-1.4543, 
 2.5134) (-1.1723, -2.6442) (-0.88892, 
 2.7435) (-0.58489, -2.8181) (-0.28622, 
 2.8614) (0.027396, -2.8766) (0.32913, 
 2.8622) (0.63949, -2.8173) (0.93194, 
 2.7457) (1.2264, -2.6426) (1.4975, 
 2.5168) (1.7640, -2.3597) (2.0028, 
 2.1848) (2.2303, -1.9802) (2.6064, -1.5196) (2.7529, 
 1.2696) (2.8767, -0.99665) (2.9671, 
 0.72392) (3.0303, -0.43292) (3.0608, 0.14854) (-3.0515, 
 0.12102) (-3.0143, -0.46756) (-2.9628, 
 0.70465) (-2.8526, -1.0370) (-2.7517, 
 1.2595) (-2.5731, -1.5641) (-2.4270, 
 1.7628) (-2.1873, -2.0271) (-2.0019, 
 2.1939) (-1.7109, -2.4073) (-1.4939, 
 2.5354) (-1.1635, -2.6889) (-0.92377, 
 2.7730) (-0.56759, -2.8605) (-0.31490, 
 2.8972) (0.052511, -2.9149) (0.30780, 
 2.9027) (0.67140, -2.8501) (0.91883, 
 2.7894) (1.2637, -2.6685) (1.4932, 
 2.5619) (1.8053, -2.3776) (2.2739, -1.9894) (2.4404, 
 1.8057) (2.6505, -1.5196) (2.7745, 
 1.3079) (2.9195, -0.98765) (2.9960, 
 0.75658) (3.0700, -0.41520) (3.0958, 0.17426) (-3.0858, 
 0.11243) (-3.0460, -0.48324) (-2.9982, 
 0.70350) (-2.8804, -1.0591) (-2.7868, 
 1.2658) (-2.5959, -1.5917) (-2.4604, 
 1.7763) (-2.2042, -2.0592) (-2.0322, 
 2.2142) (-1.7212, -2.4424) (-1.5198, 
 2.5615) (-1.1668, -2.7257) (-0.94432, 
 2.8039) (-0.56368, -2.8974) (-0.32921, 
 2.9316) (0.063444, -2.9505) (0.30032, 
 2.9392) (0.68890, -2.8829) (0.91849, 
 2.8265) (1.2871, -2.6971) (1.8336, -2.4009) (2.0211, 
 2.2632) (2.3059, -2.0064) (2.4604, 
 1.8357) (2.6849, -1.5296) (2.7999, 
 1.3329) (2.9549, -0.99025) (3.0259, 
 0.77558) (3.1050, -0.41030) (3.1290, 0.18646) (-3.1123, 
 0.11440) (-3.0724, -0.48696) (-3.0238, 
 0.71108) (-2.9056, -1.0684) (-2.8105, 
 1.2787) (-2.6189, -1.6061) (-2.4811, 
 1.7940) (-2.2239, -2.0782) (-2.0491, 
 2.2359) (-1.7369, -2.4652) (-1.5322, 
 2.5863) (-1.1779, -2.7513) (-0.95172, 
 2.8309) (-0.56965, -2.9249) (-0.33130, 
 2.9596) (0.062822, -2.9787) (0.30362, 
 2.9671) (0.69366, -2.9105) (1.2970, -2.7231) (1.5135, 
 2.6224) (1.8483, -2.4243) (2.0389, 
 2.2842) (2.3248, -2.0261) (2.4818, 
 1.8524) (2.7072, -1.5450) (2.8242, 
 1.3448) (2.9797, -1.0005) (3.0519, 
 0.78209) (3.1313, -0.41509) (3.1557, 0.18734) (-3.1284, 
 0.12317) (-3.0899, -0.48177) (-3.0377, 
 0.72308) (-2.9240, -1.0670) (-2.8216, 
 1.2934) (-2.6373, -1.6086) (-2.4890, 
 1.8108) (-2.2416, -2.0844) (-2.0535, 
 2.2542) (-1.7532, -2.4749) (-1.5329, 
 2.6053) (-1.1920, -2.7641) (-0.94869, 
 2.8498) (-0.58115, -2.9402) (-0.32467, 
 2.9776) (0.054473, -2.9960) (0.68880, -2.9291) (0.93994, 
 2.8673) (1.2959, -2.7422) (1.5288, 
 2.6337) (1.8508, -2.4430) (2.0560, 
 2.2922) (2.3310, -2.0438) (2.5000, 
 1.8568) (2.7168, -1.5608) (2.8426, 
 1.3453) (2.9923, -1.0139) (3.0699, 
 0.77878) (3.1463, -0.42553) (3.1725, 0.18032) (-3.1324, 
 0.13442) (-3.0963, -0.47161) (-3.0392, 
 0.73496) (-2.9323, -1.0584) (-2.8206, 
 1.3054) (-2.6473, -1.6018) (-2.4855, 
 1.8225) (-2.2529, -2.0798) (-2.0477, 
 2.2650) (-1.7654, -2.4726) (-1.5250, 
 2.6149) (-1.2045, -2.7643) (-0.93900, 
 2.8577) (-0.59345, -2.9428) (0.042864, -3.0008) (0.32560, 
 2.9873) (0.67836, -2.9360) (0.95239, 
 2.8687) (1.2870, -2.7510) (1.5411, 
 2.6326) (1.8440, -2.4533) (2.0678, 
 2.2887) (2.3264, -2.0551) (2.5107, 
 1.8510) (2.7145, -1.5727) (2.8519, 
 1.3376) (2.9926, -1.0260) (3.0773, 
 0.76932) (3.1491, -0.43716) (3.1777, 0.16955) (-3.1245, 
 0.14370) (-3.0905, -0.46086) (-3.0295, 
 0.74238) (-2.9289, -1.0466) (-2.8095, 
 1.3107) (-2.6464, -1.5895) (-2.4734, 
 1.8254) (-2.2546, -2.0674) (-2.0351, 
 2.2654) (-1.7695, -2.4607) (-1.5125, 
 2.6128) (-1.2110, -2.7533) (-0.60194, -2.9332) (-0.30269, 
 2.9768) (0.032673, -2.9930) (0.33499, 
 2.9785) (0.66688, -2.9303) (0.95990, 
 2.8583) (1.2747, -2.7476) (1.5465, 
 2.6211) (1.8314, -2.4524) (2.0707, 
 2.2765) (2.3140, -2.0568) (2.5111, 
 1.8386) (2.7029, -1.5768) (2.8497, 
 1.3254) (2.9821, -1.0323) (3.0727, 
 0.75799) (3.1402, -0.44556) (3.1708, 0.15948) (-3.1062, 
 0.14679) (-3.0733, -0.45396) (-3.0109, 
 0.74153) (-2.9135, -1.0361) (-2.7914, 
 1.3059) (-2.6334, -1.5759) (-2.4565, 
 1.8169) (-2.2444, -2.0513) (-2.0201, 
 2.2536) (-1.7627, -2.4427) (-1.2078, -2.7341) (-0.91755, 
 2.8362) (-0.60248, -2.9137) (-0.29661, 
 2.9582) (0.028403, -2.9740) (0.33742, 
 2.9592) (0.65906, -2.9125) (0.95857, 
 2.8389) (1.2637, -2.7317) (1.5414, 
 2.6024) (1.8175, -2.4391) (2.0621, 
 2.2594) (2.2979, -2.0466) (2.4994, 
 1.8237) (2.6852, -1.5702) (2.8353, 
 1.3134) (2.9636, -1.0295) (3.0562, 
 0.74928) (3.1217, -0.44669) (3.1530, 0.15446) (-3.0806, 
 0.14016) (-3.0468, -0.45501) (-2.9872, 
 0.72959) (-2.8872, -1.0316) (-2.7705, 
 1.2892) (-2.6084, -1.5659) (-2.4394, 
 1.7960) (-2.2218, -2.0362) (-1.7433, -2.4232) (-1.4923, 
 2.5716) (-1.1924, -2.7111) (-0.91507, 
 2.8085) (-0.59182, -2.8880) (-0.29949, 
 2.9305) (0.033954, -2.9466) (0.32929, 
 2.9325) (0.65927, -2.8846) (0.94552, 
 2.8144) (1.2585, -2.7045) (1.5240, 
 2.5811) (1.8072, -2.4137) (2.0410, 
 2.2420) (2.2829, -2.0239) (2.4755, 
 1.8112) (2.6661, -1.5513) (2.8096, 
 1.3061) (2.9412, -1.0151) (3.0297, 
 0.74752) (3.0969, -0.43736) (3.1268, 0.15832) (-3.0519, 
 0.12146) (-3.0148, -0.46721) (-2.9631, 
 0.70516) (-2.8531, -1.0368) (-2.7519, 
 1.2600) (-2.5737, -1.5639) (-2.1879, -2.0271) (-2.0018, 
 2.1946) (-1.7115, -2.4074) (-1.4937, 
 2.5360) (-1.1641, -2.6891) (-0.92348, 
 2.7736) (-0.56811, -2.8608) (-0.31450, 
 2.8977) (0.052071, -2.9154) (0.30828, 
 2.9031) (0.67105, -2.8506) (0.91938, 
 2.7897) (1.2635, -2.6691) (1.4938, 
 2.5621) (1.8052, -2.3782) (2.0080, 
 2.2295) (2.2739, -1.9900) (2.4410, 
 1.8056) (2.6506, -1.5202) (2.7750, 
 1.3077) (2.9197, -0.98819) (2.9965, 
 0.75629) (3.0703, -0.41568) (3.0963, 0.17387) (-3.0255, 
 0.090170) (-2.9823, -0.49221) (-2.9439, 
 0.66887) (-2.8159, -1.0545) (-2.5332, -1.5736) (-2.4242, 
 1.7216) (-2.1458, -2.0284) (-2.0075, 
 2.1527) (-1.6696, -2.4001) (-1.5077, 
 2.4956) (-1.1240, -2.6737) (-0.94519, 
 2.7364) (-0.53154, -2.8378) (-0.34306, 
 2.8652) (0.083641, -2.8858) (0.27405, 
 2.8767) (0.69633, -2.8156) (0.88088, 
 2.7704) (1.2814, -2.6301) (1.4526, 
 2.5507) (1.8151, -2.3369) (1.9658, 
 2.2265) (2.2753, -1.9480) (2.3995, 
 1.8111) (2.6435, -1.4793) (2.7360, 
 1.3215) (2.9044, -0.94994) (2.9614, 
 0.77777) (3.0474, -0.38170) (3.0667, 0.20218) (-3.0068, 
 0.048194) (-2.9550, -0.52932) (-2.7812, -1.0852) (-2.7409, 
 1.1738) (-2.4925, -1.5967) (-2.4342, 
 1.6758) (-2.1007, -2.0428) (-2.0268, 
 2.1093) (-1.6220, -2.4054) (-1.5354, 
 2.4565) (-1.0759, -2.6696) (-0.98027, 
 2.7031) (-0.48486, -2.8245) (-0.38405, 
 2.8392) (0.12698, -2.8638) (0.22883, 
 2.8589) (0.73456, -2.7858) (0.83327, 
 2.7617) (1.3130, -2.5938) (1.4045, 
 2.5513) (1.8387, -2.2955) (1.9193, 
 2.2364) (2.2900, -1.9032) (2.3564, 
 1.8300) (2.6486, -1.4329) (2.6981, 
 1.3486) (2.8998, -0.90390) (2.9303, 
 0.81188) (3.0332, -0.33790) (3.0435, 
 0.24195) (-2.9381, -0.57507) (-2.9381, 
 0.57507) (-2.7548, -1.1266) (-2.7548, 
 1.1266) (-2.4577, -1.6321) (-2.4577, 
 1.6321) (-2.0591, -2.0708) (-2.0591, 
 2.0708) (-1.5751, -2.4248) (-1.5751, 
 2.4248) (-1.0257, -2.6796) (-1.0257, 
 2.6796) (-0.43342, -2.8247) (-0.43342, 
 2.8247) (0.17756, -2.8542) (0.17756, 
 2.8542) (0.78221, -2.7667) (0.78221, 
 2.7667) (1.3558, -2.5660) (1.3558, 
 2.5660) (1.8748, -2.2602) (1.8748, 
 2.2602) (2.3181, -1.8619) (2.3181, 
 1.8619) (2.6674, -1.3872) (2.6674, 
 1.3872) (2.9085, -0.85573) (2.9085, 
 0.85573) (3.0316, -0.28919) (3.0316, 
 0.28919) (-3.0068, -0.048194) (-2.9550, 0.52932) (-2.7812, 
 1.0852) (-2.7409, -1.1738) (-2.4925, 
 1.5967) (-2.4342, -1.6758) (-2.1007, 
 2.0428) (-2.0268, -2.1093) (-1.6220, 
 2.4054) (-1.5354, -2.4565) (-1.0759, 
 2.6696) (-0.98027, -2.7031) (-0.48486, 
 2.8245) (-0.38405, -2.8392) (0.12698, 
 2.8638) (0.22883, -2.8589) (0.73456, 
 2.7858) (0.83327, -2.7617) (1.3130, 
 2.5938) (1.4045, -2.5513) (1.8387, 
 2.2955) (1.9193, -2.2364) (2.2900, 
 1.9032) (2.3564, -1.8300) (2.6486, 
 1.4329) (2.6981, -1.3486) (2.8998, 
 0.90390) (2.9303, -0.81188) (3.0332, 
 0.33790) (3.0435, -0.24195) (-3.0255, -0.090170) (-2.9823, 
 0.49221) (-2.9439, -0.66887) (-2.8159, 1.0545) (-2.5332, 
 1.5736) (-2.4242, -1.7216) (-2.1458, 
 2.0284) (-2.0075, -2.1527) (-1.6696, 
 2.4001) (-1.5077, -2.4956) (-1.1240, 
 2.6737) (-0.94519, -2.7364) (-0.53154, 
 2.8378) (-0.34306, -2.8652) (0.083641, 
 2.8858) (0.27405, -2.8767) (0.69633, 
 2.8156) (0.88088, -2.7704) (1.2814, 
 2.6301) (1.4526, -2.5507) (1.8151, 
 2.3369) (1.9658, -2.2265) (2.2753, 
 1.9480) (2.3995, -1.8111) (2.6435, 
 1.4793) (2.7360, -1.3215) (2.9044, 
 0.94994) (2.9614, -0.77777) (3.0474, 
 0.38170) (3.0667, -0.20218) (-3.0519, -0.12146) (-3.0148, 
 0.46721) (-2.9631, -0.70516) (-2.8531, 
 1.0368) (-2.7519, -1.2600) (-2.5737, 1.5639) (-2.1879, 
 2.0271) (-2.0018, -2.1946) (-1.7115, 
 2.4074) (-1.4937, -2.5360) (-1.1641, 
 2.6891) (-0.92348, -2.7736) (-0.56811, 
 2.8608) (-0.31450, -2.8977) (0.052071, 
 2.9154) (0.30828, -2.9031) (0.67105, 
 2.8506) (0.91938, -2.7897) (1.2635, 
 2.6691) (1.4938, -2.5621) (1.8052, 
 2.3782) (2.0080, -2.2295) (2.2739, 
 1.9900) (2.4410, -1.8056) (2.6506, 
 1.5202) (2.7750, -1.3077) (2.9197, 
 0.98819) (2.9965, -0.75629) (3.0703, 
 0.41568) (3.0963, -0.17387) (-3.0806, -0.14016) (-3.0468, 
 0.45501) (-2.9872, -0.72959) (-2.8872, 
 1.0316) (-2.7705, -1.2892) (-2.6084, 
 1.5659) (-2.4394, -1.7960) (-2.2218, 2.0362) (-1.7433, 
 2.4232) (-1.4923, -2.5716) (-1.1924, 
 2.7111) (-0.91507, -2.8085) (-0.59182, 
 2.8880) (-0.29949, -2.9305) (0.033954, 
 2.9466) (0.32929, -2.9325) (0.65927, 
 2.8846) (0.94552, -2.8144) (1.2585, 
 2.7045) (1.5240, -2.5811) (1.8072, 
 2.4137) (2.0410, -2.2420) (2.2829, 
 2.0239) (2.4755, -1.8112) (2.6661, 
 1.5513) (2.8096, -1.3061) (2.9412, 
 1.0151) (3.0297, -0.74752) (3.0969, 
 0.43736) (3.1268, -0.15832) (-3.1062, -0.14679) (-3.0733, 
 0.45396) (-3.0109, -0.74153) (-2.9135, 
 1.0361) (-2.7914, -1.3059) (-2.6334, 
 1.5759) (-2.4565, -1.8169) (-2.2444, 
 2.0513) (-2.0201, -2.2536) (-1.7627, 2.4427) (-1.2078, 
 2.7341) (-0.91755, -2.8362) (-0.60248, 
 2.9137) (-0.29661, -2.9582) (0.028403, 
 2.9740) (0.33742, -2.9592) (0.65906, 
 2.9125) (0.95857, -2.8389) (1.2637, 
 2.7317) (1.5414, -2.6024) (1.8175, 
 2.4391) (2.0621, -2.2594) (2.2979, 
 2.0466) (2.4994, -1.8237) (2.6852, 
 1.5702) (2.8353, -1.3134) (2.9636, 
 1.0295) (3.0562, -0.74928) (3.1217, 
 0.44669) (3.1530, -0.15446) (-3.1245, -0.14370) (-3.0905, 
 0.46086) (-3.0295, -0.74238) (-2.9289, 
 1.0466) (-2.8095, -1.3107) (-2.6464, 
 1.5895) (-2.4734, -1.8254) (-2.2546, 
 2.0674) (-2.0351, -2.2654) (-1.7695, 
 2.4607) (-1.5125, -2.6128) (-1.2110, 2.7533) (-0.60194, 
 2.9332) (-0.30269, -2.9768) (0.032673, 
 2.9930) (0.33499, -2.9785) (0.66688, 
 2.9303) (0.95990, -2.8583) (1.2747, 
 2.7476) (1.5465, -2.6211) (1.8314, 
 2.4524) (2.0707, -2.2765) (2.3140, 
 2.0568) (2.5111, -1.8386) (2.7029, 
 1.5768) (2.8497, -1.3254) (2.9821, 
 1.0323) (3.0727, -0.75799) (3.1402, 
 0.44556) (3.1708, -0.15948) (-3.1324, -0.13442) (-3.0963, 
 0.47161) (-3.0392, -0.73496) (-2.9323, 
 1.0584) (-2.8206, -1.3054) (-2.6473, 
 1.6018) (-2.4855, -1.8225) (-2.2529, 
 2.0798) (-2.0477, -2.2650) (-1.7654, 
 2.4726) (-1.5250, -2.6149) (-1.2045, 
 2.7643) (-0.93900, -2.8577) (-0.59345, 2.9428) (0.042864, 
 3.0008) (0.32560, -2.9873) (0.67836, 
 2.9360) (0.95239, -2.8687) (1.2870, 
 2.7510) (1.5411, -2.6326) (1.8440, 
 2.4533) (2.0678, -2.2887) (2.3264, 
 2.0551) (2.5107, -1.8510) (2.7145, 
 1.5727) (2.8519, -1.3376) (2.9926, 
 1.0260) (3.0773, -0.76932) (3.1491, 
 0.43716) (3.1777, -0.16955) (-3.1284, -0.12317) (-3.0899, 
 0.48177) (-3.0377, -0.72308) (-2.9240, 
 1.0670) (-2.8216, -1.2934) (-2.6373, 
 1.6086) (-2.4890, -1.8108) (-2.2416, 
 2.0844) (-2.0535, -2.2542) (-1.7532, 
 2.4749) (-1.5329, -2.6053) (-1.1920, 
 2.7641) (-0.94869, -2.8498) (-0.58115, 
 2.9402) (-0.32467, -2.9776) (0.054473, 2.9960) (0.68880, 
 2.9291) (0.93994, -2.8673) (1.2959, 
 2.7422) (1.5288, -2.6337) (1.8508, 
 2.4430) (2.0560, -2.2922) (2.3310, 
 2.0438) (2.5000, -1.8568) (2.7168, 
 1.5608) (2.8426, -1.3453) (2.9923, 
 1.0139) (3.0699, -0.77878) (3.1463, 
 0.42553) (3.1725, -0.18032) (-3.1123, -0.11440) (-3.0724, 
 0.48696) (-3.0238, -0.71108) (-2.9056, 
 1.0684) (-2.8105, -1.2787) (-2.6189, 
 1.6061) (-2.4811, -1.7940) (-2.2239, 
 2.0782) (-2.0491, -2.2359) (-1.7369, 
 2.4652) (-1.5322, -2.5863) (-1.1779, 
 2.7513) (-0.95172, -2.8309) (-0.56965, 
 2.9249) (-0.33130, -2.9596) (0.062822, 
 2.9787) (0.30362, -2.9671) (0.69366, 2.9105) (1.2970, 
 2.7231) (1.5135, -2.6224) (1.8483, 
 2.4243) (2.0389, -2.2842) (2.3248, 
 2.0261) (2.4818, -1.8524) (2.7072, 
 1.5450) (2.8242, -1.3448) (2.9797, 
 1.0005) (3.0519, -0.78209) (3.1313, 
 0.41509) (3.1557, -0.18734) (-3.0858, -0.11243) (-3.0460, 
 0.48324) (-2.9982, -0.70350) (-2.8804, 
 1.0591) (-2.7868, -1.2658) (-2.5959, 
 1.5917) (-2.4604, -1.7763) (-2.2042, 
 2.0592) (-2.0322, -2.2142) (-1.7212, 
 2.4424) (-1.5198, -2.5615) (-1.1668, 
 2.7257) (-0.94432, -2.8039) (-0.56368, 
 2.8974) (-0.32921, -2.9316) (0.063444, 
 2.9505) (0.30032, -2.9392) (0.68890, 
 2.8829) (0.91849, -2.8265) (1.2871, 2.6971) (1.8336, 
 2.4009) (2.0211, -2.2632) (2.3059, 
 2.0064) (2.4604, -1.8357) (2.6849, 
 1.5296) (2.7999, -1.3329) (2.9549, 
 0.99025) (3.0259, -0.77558) (3.1050, 
 0.41030) (3.1290, -0.18646) (-3.0515, -0.12102) (-3.0143, 
 0.46756) (-2.9628, -0.70465) (-2.8526, 
 1.0370) (-2.7517, -1.2595) (-2.5731, 
 1.5641) (-2.4270, -1.7628) (-2.1873, 
 2.0271) (-2.0019, -2.1939) (-1.7109, 
 2.4073) (-1.4939, -2.5354) (-1.1635, 
 2.6889) (-0.92377, -2.7730) (-0.56759, 
 2.8605) (-0.31490, -2.8972) (0.052511, 
 2.9149) (0.30780, -2.9027) (0.67140, 
 2.8501) (0.91883, -2.7894) (1.2637, 
 2.6685) (1.4932, -2.5619) (1.8053, 2.3776) (2.2739, 
 1.9894) (2.4404, -1.8057) (2.6505, 
 1.5196) (2.7745, -1.3079) (2.9195, 
 0.98765) (2.9960, -0.75658) (3.0700, 
 0.41520) (3.0958, -0.17426) (-3.0140, -0.14284) (-2.9822, 
 0.43822) (-2.9214, -0.71806) (-2.8273, 
 1.0014) (-2.7081, -1.2639) (-2.5556, 
 1.5236) (-2.3829, -1.7581) (-2.1783, 
 1.9834) (-1.9592, -2.1803) (-1.7108, 
 2.3622) (-1.4543, -2.5134) (-1.1723, 
 2.6442) (-0.88892, -2.7435) (-0.58489, 
 2.8181) (-0.28622, -2.8614) (0.027396, 
 2.8766) (0.32913, -2.8622) (0.63949, 
 2.8173) (0.93194, -2.7457) (1.2264, 
 2.6426) (1.4975, -2.5168) (1.7640, 
 2.3597) (2.0028, -2.1848) (2.2303, 1.9802) (2.6064, 
 1.5196) (2.7529, -1.2696) (2.8767, 
 0.99665) (2.9671, -0.72392) (3.0303, 
 0.43292) (3.0608, -0.14854) (-2.9787, -0.17879) (-2.9554, 
 0.39554) (-2.8791, -0.74581) (-2.8100, 
 0.95370) (-2.6605, -1.2823) (-2.5485, 
 1.4729) (-2.3321, -1.7664) (-2.1817, 
 1.9318) (-1.9071, -2.1782) (-1.7246, 
 2.3117) (-1.4031, -2.5009) (-1.1959, 
 2.5970) (-0.84076, -2.7213) (-0.61737, 
 2.7760) (-0.24303, -2.8302) (-0.012599, 
 2.8413) (0.36557, -2.8233) (0.59362, 
 2.7904) (0.96015, -2.7008) (1.1765, 
 2.6252) (1.5164, -2.4677) (1.7121, 
 2.3524) (2.0115, -2.1335) (2.1786, 
 1.9834) (2.4252, -1.7119) (2.5569, 1.5330) (2.8316, 
 1.0199) (2.9449, -0.67851) (2.9913, 
 0.46496) (3.0296, -0.10900) (-2.9516, -0.22729) (-2.9392, 
 0.34234) (-2.8421, -0.78762) (-2.8054, 
 0.89797) (-2.6153, -1.3157) (-2.5557, 
 1.4169) (-2.2803, -1.7901) (-2.2004, 
 1.8778) (-1.8510, -2.1911) (-1.7540, 
 2.2620) (-1.3450, -2.5026) (-1.2349, 
 2.5536) (-0.78295, -2.7116) (-0.66426, 
 2.7406) (-0.18792, -2.8096) (-0.065490, 
 2.8155) (0.41572, -2.7926) (0.53689, 
 2.7752) (1.0033, -2.6613) (1.1182, 
 2.6211) (1.5507, -2.4210) (1.6547, 
 2.3598) (2.0356, -2.0815) (2.1245, 
 2.0017) (2.4382, -1.6567) (2.5082, 
 1.5617) (2.7419, -1.1641) (2.7903, 1.0577) (2.9591, 
 0.51039) (3.0078, -0.057860)};
\draw[very thick] (0,0) circle (50mm);
\draw[-stealth, dashed,gray] (0,0) -- (32:50mm) 
node[below left, rotate=30, black]{$2+|\xi|$}; 
\foreach \x in {-4,...,-1}
\draw[shift={(\x,0)}] (0pt, 1pt) -- (0pt,-1pt) node[below] {\tiny $\x$};
\foreach \x in {1,...,4}
\draw[shift={(\x,0)}] (0pt, 1pt) -- (0pt,-1pt) node[below] {\tiny $\x$};
\end{tikzpicture}
\end{minipage}
\hspace{3.5mm}
\begin{minipage}{0.46\textwidth}
\begin{tikzpicture}[domain=-0.7:0.7,font=\sffamily, scale=0.55]
\draw[-stealth] (-5.3,0) -- (5.3,0) node[above] {$x$};
\draw[-stealth] (0,-5.3) -- (0,5.3) node[left] {$y$};
\draw plot [only marks, mark=*, mark options={fill=black},mark size=0.3pt]
coordinates{
(-2.9420, -0.28485) (-2.9420, 
 0.28485) (-2.8206, -0.84290) (-2.8206, 
 0.84290) (-2.5828, -1.3665) (-2.5828, 
 1.3665) (-2.2383, -1.8342) (-2.2383, 
 1.8342) (-1.8013, -2.2269) (-1.8013, 
 2.2269) (-1.2897, -2.5284) (-1.2897, 
 2.5284) (-0.72449, -2.7265) (-0.72449, 
 2.7265) (-0.12876, -2.8130) (-0.12876, 
 2.8130) (0.47309, -2.7843) (0.47309, 
 2.7843) (1.0564, -2.6415) (1.0564, 
 2.6415) (1.5974, -2.3906) (1.5974, 
 2.3906) (2.0738, -2.0418) (2.0738, 
 2.0418) (2.4663, -1.6093) (2.4663, 
 1.6093) (2.7587, -1.1109) (2.7587, 
 1.1109) (2.9390, -0.56708) (2.9390, 0.56708) (-2.9706, 
 0.039968) (-2.9178, -0.53046) (-2.9008, 
 0.60876) (-2.7446, -1.0792) (-2.7113, 
 1.1527) (-2.4581, -1.5838) (-2.4098, 
 1.6494) (-2.0701, -2.0237) (-2.0088, 
 2.0788) (-1.5965, -2.3807) (-1.5247, 
 2.4230) (-1.0567, -2.6403) (-0.97739, 
 2.6681) (-0.47271, -2.7919) (-0.38917, 
 2.8040) (0.13142, -2.8291) (0.21582, 
 2.8251) (0.73102, -2.7505) (0.81282, 
 2.7305) (1.3016, -2.5593) (1.3774, 
 2.5240) (1.8197, -2.2632) (1.8865, 
 2.2142) (2.2642, -1.8744) (2.3192, 
 1.8138) (2.6169, -1.4089) (2.6580, 
 1.3390) (2.8635, -0.88570) (2.9937, -0.32621) (3.0022, 
 0.24664) (-2.9993, -0.19399) (-2.9787, 
 0.38587) (-2.8962, -0.76591) (-2.8352, 
 0.94994) (-2.6736, -1.3065) (-2.5747, 
 1.4752) (-2.3405, -1.7937) (-2.2078, 
 1.9401) (-1.9108, -2.2075) (-1.7497, 
 2.3256) (-1.4020, -2.5310) (-1.2191, 
 2.6160) (-0.83493, -2.7509) (-0.63778, 
 2.7993) (-0.23285, -2.8581) (-0.029479, 
 2.8680) (0.37961, -2.8484) (0.58087, 
 2.8192) (0.97738, -2.7220) (1.1683, 
 2.6550) (1.5360, -2.4842) (1.7088, 
 2.3821) (2.0326, -2.1446) (2.1802, 
 2.0116) (2.4470, -1.7171) (2.7620, -1.2194) (2.8423, 
 1.0420) (2.9649, -0.67168) (3.0059, 
 0.48266) (3.0473, -0.096483) (-3.0455, 
 0.17840) (-3.0206, -0.41092) (-2.9448, 
 0.76043) (-2.8710, -0.98344) (-2.7226, 
 1.3114) (-2.6030, -1.5157) (-2.3880, 
 1.8087) (-2.2275, -1.9860) (-1.9548, 
 2.2320) (-1.7599, -2.3751) (-1.4407, 
 2.5640) (-1.2195, -2.6670) (-0.86672, 
 2.7910) (-0.62821, -2.8497) (-0.25644, 
 2.9038) (-0.010400, -2.9157) (0.36519, 
 2.8977) (0.60868, -2.8623) (0.97272, 
 2.7729) (1.2037, -2.6918) (1.5413, 
 2.5346) (1.7503, -2.4110) (2.2262, -2.0314) (2.4711, 
 1.7606) (2.6118, -1.5687) (2.7943, 
 1.2566) (2.8914, -1.0417) (3.0039, 
 0.70110) (3.0535, -0.47206) (3.0916, 
 0.11692) (-3.0929, -0.013962) (-3.0321, 
 0.58194) (-3.0262, -0.60929) (-2.8462, 
 1.1540) (-2.8346, -1.1797) (-2.5429, 
 1.6789) (-2.5261, -1.7019) (-2.1345, 
 2.1352) (-2.1132, -2.1544) (-1.6379, 
 2.5040) (-1.6129, -2.5188) (-1.0733, 
 2.7704) (-1.0458, -2.7801) (-0.46404, 
 2.9234) (-0.43499, -2.9276) (0.16511, 
 2.9567) (0.19446, -2.9553) (0.78836, 
 2.8689) (0.81681, -2.8619) (1.4066, -2.6514) (1.9164, 
 2.3493) (1.9396, -2.3322) (2.3751, 
 1.9387) (2.3942, -1.9176) (2.7374, 
 1.4488) (2.7516, -1.4244) (2.9885, 
 0.89950) (2.9973, -0.87285) (3.1183, 
 0.31338) (3.1212, -0.28558) (-3.1150, -0.19430) (-3.0920, 
 0.41011) (-3.0097, -0.79075) (-2.9413, 
 0.99774) (-2.7801, -1.3549) (-2.6694, 
 1.5446) (-2.4358, -1.8636) (-2.2872, 
 2.0282) (-1.9909, -2.2960) (-1.8105, 
 2.4289) (-1.4637, -2.6346) (-1.2589, 
 2.7302) (-0.87565, -2.8652) (-0.65482, 
 2.9197) (-0.25092, -2.9786) (-0.023123, 
 2.9897) (0.38493, -2.9701) (1.0059, -2.8399) (1.2197, 
 2.7645) (1.5865, -2.5934) (1.7800, 
 2.4786) (2.1031, -2.2407) (2.2684, 
 2.0912) (2.5344, -1.7962) (2.6647, 
 1.6181) (2.8630, -1.2782) (2.9529, 
 1.0788) (3.0752, -0.70786) (3.1211, 
 0.49528) (3.1624, -0.10851) (-3.1219, 
 0.23983) (-3.1083, -0.36703) (-3.0067, 
 0.83688) (-2.9665, -0.95888) (-2.7676, 
 1.3997) (-2.7023, -1.5115) (-2.4142, 
 1.9053) (-2.3266, -2.0023) (-1.9610, 
 2.3329) (-1.8547, -2.4112) (-1.4268, 
 2.6651) (-1.3061, -2.7214) (-0.83322, 
 2.8882) (-0.70310, -2.9202) (-0.070502, -2.9995) (0.43299, 
 2.9753) (0.56583, -2.9560) (1.0538, 
 2.8358) (1.1799, -2.7914) (1.6324, 
 2.5801) (1.7465, -2.5125) (2.1451, 
 2.2188) (2.2424, -2.1307) (2.5708, 
 1.7666) (2.6475, -1.6616) (2.8921, 
 1.2421) (2.9451, -1.1245) (3.0960, 
 0.66666) (3.1231, -0.54136) (3.1741, 0.063960) (-3.1219, 
 0.065897) (-3.0714, -0.53721) (-3.0435, 
 0.66630) (-2.8943, -1.1183) (-2.8397, 
 1.2395) (-2.5978, -1.6537) (-2.5186, 
 1.7619) (-2.1940, -2.1215) (-2.0935, 
 2.2123) (-1.6995, -2.5024) (-1.5819, 
 2.5722) (-1.1346, -2.7810) (-0.52239, -2.9457) (-0.38541, 
 2.9657) (0.11205, -2.9898) (0.25044, 
 2.9831) (0.74276, -2.9115) (0.87690, 
 2.8784) (1.3439, -2.7139) (1.4683, 
 2.6559) (1.8910, -2.4052) (2.0005, 
 2.3245) (2.3615, -1.9980) (2.4517, 
 1.8980) (2.7362, -1.5090) (2.8034, 
 1.3937) (2.9998, -0.95820) (3.0413, 
 0.83240) (3.1415, -0.36815) (3.1555, 
 0.23698) (-3.0940, -0.11651) (-3.0549, 
 0.48176) (-3.0055, -0.71001) (-2.8898, 
 1.0603) (-2.7930, -1.2745) (-2.6054, 
 1.5955) (-2.4652, -1.7868) (-2.2135, 
 2.0655) (-2.0357, -2.2261) (-1.5219, -2.5742) (-1.1749, 
 2.7360) (-0.94494, -2.8170) (-0.57087, 
 2.9092) (-0.32847, -2.9445) (0.057396, 
 2.9632) (0.30229, -2.9515) (0.68417, 
 2.8959) (0.92153, -2.8375) (1.2838, 
 2.7100) (1.5039, -2.6074) (1.8317, 
 2.4132) (2.0256, -2.2705) (2.3056, 
 2.0175) (2.4653, -1.8406) (2.6859, 
 1.5391) (2.8049, -1.3353) (2.9573, 
 0.99779) (3.0306, -0.77537) (3.1084, 
 0.41558) (3.1332, -0.18365) (-3.0461, 
 0.28040) (-3.0428, -0.31131) (-2.9238, 
 0.86065) (-2.9140, -0.89029) (-2.6808, 
 1.4057) (-2.6649, -1.4329) (-2.3058, -1.9168) (-1.8773, 
 2.3034) (-1.8514, -2.3224) (-1.3498, 
 2.6192) (-1.3204, -2.6329) (-0.76618, 
 2.8279) (-0.73449, -2.8357) (-0.15035, 
 2.9208) (-0.11765, -2.9223) (0.47249, 
 2.8941) (0.50485, -2.8894) (1.0768, 
 2.7488) (1.1075, -2.7381) (1.6380, 
 2.4910) (1.6658, -2.4746) (2.1330, 
 2.1312) (2.1567, -2.1098) (2.5415, 
 1.6841) (2.5602, -1.6586) (2.8469, 
 1.1681) (2.8598, -1.1395) (3.0367, 
 0.60418) (3.0433, -0.57374) (3.1030, 0.015541) (-3.0244, 
 0.061831) (-2.9751, -0.52050) (-2.9488, 
 0.64163) (-2.8031, -1.0816) (-2.5154, -1.5984) (-2.4409, 
 1.6999) (-2.1239, -2.0498) (-2.0293, 
 2.1351) (-1.6445, -2.4174) (-1.5337, 
 2.4829) (-1.0969, -2.6861) (-0.97454, 
 2.7290) (-0.50359, -2.8448) (-0.37462, 
 2.8635) (0.11118, -2.8870) (0.24148, 
 2.8808) (0.72225, -2.8110) (0.84854, 
 2.7800) (1.3046, -2.6199) (1.4217, 
 2.5655) (1.8345, -2.3215) (1.9376, 
 2.2458) (2.2901, -1.9281) (2.3750, 
 1.8342) (2.6529, -1.4556) (2.7162, 
 1.3475) (2.9079, -0.92359) (2.9470, 
 0.80554) (3.0449, -0.35372) (3.0581, 
 0.23063) (-2.9976, -0.17048) (-2.8995, -0.74248) (-2.8239, 
 0.97077) (-2.6818, -1.2841) (-2.5591, 
 1.4934) (-2.3533, -1.7733) (-2.1887, 
 1.9549) (-1.9276, -2.1899) (-1.7278, 
 2.3364) (-1.4221, -2.5169) (-1.1953, 
 2.6223) (-0.85760, -2.7408) (-0.61300, 
 2.8009) (-0.25711, -2.8526) (-0.0047952, 
 2.8648) (0.35475, -2.8476) (0.60446, 
 2.8114) (0.95295, -2.7260) (1.1898, 
 2.6429) (1.5130, -2.4927) (1.7274, 
 2.3662) (2.0120, -2.1574) (2.1951, 
 1.9925) (2.4295, -1.7337) (2.5738, 
 1.5372) (2.7485, -1.2390) (2.8481, 
 1.0190) (2.9558, -0.69353) (3.0067, 
 0.45906) (3.0431, -0.11969) (-2.9976, 0.17048) (-2.8995, 
 0.74248) (-2.8239, -0.97077) (-2.6818, 
 1.2841) (-2.5591, -1.4934) (-2.3533, 
 1.7733) (-2.1887, -1.9549) (-1.9276, 
 2.1899) (-1.7278, -2.3364) (-1.4221, 
 2.5169) (-1.1953, -2.6223) (-0.85760, 
 2.7408) (-0.61300, -2.8009) (-0.25711, 
 2.8526) (-0.0047952, -2.8648) (0.35475, 
 2.8476) (0.60446, -2.8114) (0.95295, 
 2.7260) (1.1898, -2.6429) (1.5130, 
 2.4927) (1.7274, -2.3662) (2.0120, 
 2.1574) (2.1951, -1.9925) (2.4295, 
 1.7337) (2.5738, -1.5372) (2.7485, 
 1.2390) (2.8481, -1.0190) (2.9558, 
 0.69353) (3.0067, -0.45906) (3.0431, 
 0.11969) (-3.0244, -0.061831) (-2.9751, 
 0.52050) (-2.9488, -0.64163) (-2.8031, 1.0816) (-2.5154, 
 1.5984) (-2.4409, -1.6999) (-2.1239, 
 2.0498) (-2.0293, -2.1351) (-1.6445, 
 2.4174) (-1.5337, -2.4829) (-1.0969, 
 2.6861) (-0.97454, -2.7290) (-0.50359, 
 2.8448) (-0.37462, -2.8635) (0.11118, 
 2.8870) (0.24148, -2.8808) (0.72225, 
 2.8110) (0.84854, -2.7800) (1.3046, 
 2.6199) (1.4217, -2.5655) (1.8345, 
 2.3215) (1.9376, -2.2458) (2.2901, 
 1.9281) (2.3750, -1.8342) (2.6529, 
 1.4556) (2.7162, -1.3475) (2.9079, 
 0.92359) (2.9470, -0.80554) (3.0449, 
 0.35372) (3.0581, -0.23063) (-3.0461, -0.28040) (-3.0428, 
 0.31131) (-2.9238, -0.86065) (-2.9140, 
 0.89029) (-2.6808, -1.4057) (-2.6649, 1.4329) (-2.3058, 
 1.9168) (-1.8773, -2.3034) (-1.8514, 
 2.3224) (-1.3498, -2.6192) (-1.3204, 
 2.6329) (-0.76618, -2.8279) (-0.73449, 
 2.8357) (-0.15035, -2.9208) (-0.11765, 
 2.9223) (0.47249, -2.8941) (0.50485, 
 2.8894) (1.0768, -2.7488) (1.1075, 
 2.7381) (1.6380, -2.4910) (1.6658, 
 2.4746) (2.1330, -2.1312) (2.1567, 
 2.1098) (2.5415, -1.6841) (2.5602, 
 1.6586) (2.8469, -1.1681) (2.8598, 
 1.1395) (3.0367, -0.60418) (3.0433, 
 0.57374) (3.1030, -0.015541) (-3.0940, 
 0.11651) (-3.0549, -0.48176) (-3.0055, 
 0.71001) (-2.8898, -1.0603) (-2.7930, 
 1.2745) (-2.6054, -1.5955) (-2.4652, 
 1.7868) (-2.2135, -2.0655) (-2.0357, 2.2261) (-1.5219, 
 2.5742) (-1.1749, -2.7360) (-0.94494, 
 2.8170) (-0.57087, -2.9092) (-0.32847, 
 2.9445) (0.057396, -2.9632) (0.30229, 
 2.9515) (0.68417, -2.8959) (0.92153, 
 2.8375) (1.2838, -2.7100) (1.5039, 
 2.6074) (1.8317, -2.4132) (2.0256, 
 2.2705) (2.3056, -2.0175) (2.4653, 
 1.8406) (2.6859, -1.5391) (2.8049, 
 1.3353) (2.9573, -0.99779) (3.0306, 
 0.77537) (3.1084, -0.41558) (3.1332, 
 0.18365) (-3.1219, -0.065897) (-3.0714, 
 0.53721) (-3.0435, -0.66630) (-2.8943, 
 1.1183) (-2.8397, -1.2395) (-2.5978, 
 1.6537) (-2.5186, -1.7619) (-2.1940, 
 2.1215) (-2.0935, -2.2123) (-1.6995, 
 2.5024) (-1.5819, -2.5722) (-1.1346, 2.7810) (-0.52239, 
 2.9457) (-0.38541, -2.9657) (0.11205, 
 2.9898) (0.25044, -2.9831) (0.74276, 
 2.9115) (0.87690, -2.8784) (1.3439, 
 2.7139) (1.4683, -2.6559) (1.8910, 
 2.4052) (2.0005, -2.3245) (2.3615, 
 1.9980) (2.4517, -1.8980) (2.7362, 
 1.5090) (2.8034, -1.3937) (2.9998, 
 0.95820) (3.0413, -0.83240) (3.1415, 
 0.36815) (3.1555, -0.23698) (-3.1219, -0.23983) (-3.1083, 
 0.36703) (-3.0067, -0.83688) (-2.9665, 
 0.95888) (-2.7676, -1.3997) (-2.7023, 
 1.5115) (-2.4142, -1.9053) (-2.3266, 
 2.0023) (-1.9610, -2.3329) (-1.8547, 
 2.4112) (-1.4268, -2.6651) (-1.3061, 
 2.7214) (-0.83322, -2.8882) (-0.70310, 2.9202) (-0.070502, 
 2.9995) (0.43299, -2.9753) (0.56583, 
 2.9560) (1.0538, -2.8358) (1.1799, 
 2.7914) (1.6324, -2.5801) (1.7465, 
 2.5125) (2.1451, -2.2188) (2.2424, 
 2.1307) (2.5708, -1.7666) (2.6475, 
 1.6616) (2.8921, -1.2421) (2.9451, 
 1.1245) (3.0960, -0.66666) (3.1231, 
 0.54136) (3.1741, -0.063960) (-3.1150, 
 0.19430) (-3.0920, -0.41011) (-3.0097, 
 0.79075) (-2.9413, -0.99774) (-2.7801, 
 1.3549) (-2.6694, -1.5446) (-2.4358, 
 1.8636) (-2.2872, -2.0282) (-1.9909, 
 2.2960) (-1.8105, -2.4289) (-1.4637, 
 2.6346) (-1.2589, -2.7302) (-0.87565, 
 2.8652) (-0.65482, -2.9197) (-0.25092, 
 2.9786) (-0.023123, -2.9897) (0.38493, 2.9701) (1.0059, 
 2.8399) (1.2197, -2.7645) (1.5865, 
 2.5934) (1.7800, -2.4786) (2.1031, 
 2.2407) (2.2684, -2.0912) (2.5344, 
 1.7962) (2.6647, -1.6181) (2.8630, 
 1.2782) (2.9529, -1.0788) (3.0752, 
 0.70786) (3.1211, -0.49528) (3.1624, 0.10851) (-3.0929, 
 0.013962) (-3.0321, -0.58194) (-3.0262, 
 0.60929) (-2.8462, -1.1540) (-2.8346, 
 1.1797) (-2.5429, -1.6789) (-2.5261, 
 1.7019) (-2.1345, -2.1352) (-2.1132, 
 2.1544) (-1.6379, -2.5040) (-1.6129, 
 2.5188) (-1.0733, -2.7704) (-1.0458, 
 2.7801) (-0.46404, -2.9234) (-0.43499, 
 2.9276) (0.16511, -2.9567) (0.19446, 
 2.9553) (0.78836, -2.8689) (0.81681, 2.8619) (1.4066, 
 2.6514) (1.9164, -2.3493) (1.9396, 
 2.3322) (2.3751, -1.9387) (2.3942, 
 1.9176) (2.7374, -1.4488) (2.7516, 
 1.4244) (2.9885, -0.89950) (2.9973, 
 0.87285) (3.1183, -0.31338) (3.1212, 
 0.28558) (-3.0455, -0.17840) (-3.0206, 
 0.41092) (-2.9448, -0.76043) (-2.8710, 
 0.98344) (-2.7226, -1.3114) (-2.6030, 
 1.5157) (-2.3880, -1.8087) (-2.2275, 
 1.9860) (-1.9548, -2.2320) (-1.7599, 
 2.3751) (-1.4407, -2.5640) (-1.2195, 
 2.6670) (-0.86672, -2.7910) (-0.62821, 
 2.8497) (-0.25644, -2.9038) (-0.010400, 
 2.9157) (0.36519, -2.8977) (0.60868, 
 2.8623) (0.97272, -2.7729) (1.2037, 
 2.6918) (1.5413, -2.5346) (1.7503, 2.4110) (2.2262, 
 2.0314) (2.4711, -1.7606) (2.6118, 
 1.5687) (2.7943, -1.2566) (2.8914, 
 1.0417) (3.0039, -0.70110) (3.0535, 
 0.47206) (3.0916, -0.11692) (-2.9993, 
 0.19399) (-2.9787, -0.38587) (-2.8962, 
 0.76591) (-2.8352, -0.94994) (-2.6736, 
 1.3065) (-2.5747, -1.4752) (-2.3405, 
 1.7937) (-2.2078, -1.9401) (-1.9108, 
 2.2075) (-1.7497, -2.3256) (-1.4020, 
 2.5310) (-1.2191, -2.6160) (-0.83493, 
 2.7509) (-0.63778, -2.7993) (-0.23285, 
 2.8581) (-0.029479, -2.8680) (0.37961, 
 2.8484) (0.58087, -2.8192) (0.97738, 
 2.7220) (1.1683, -2.6550) (1.5360, 
 2.4842) (1.7088, -2.3821) (2.0326, 
 2.1446) (2.1802, -2.0116) (2.4470, 1.7171) (2.7620, 
 1.2194) (2.8423, -1.0420) (2.9649, 
 0.67168) (3.0059, -0.48266) (3.0473, 
 0.096483) (-2.9706, -0.039968) (-2.9178, 
 0.53046) (-2.9008, -0.60876) (-2.7446, 
 1.0792) (-2.7113, -1.1527) (-2.4581, 
 1.5838) (-2.4098, -1.6494) (-2.0701, 
 2.0237) (-2.0088, -2.0788) (-1.5965, 
 2.3807) (-1.5247, -2.4230) (-1.0567, 
 2.6403) (-0.97739, -2.6681) (-0.47271, 
 2.7919) (-0.38917, -2.8040) (0.13142, 
 2.8291) (0.21582, -2.8251) (0.73102, 
 2.7505) (0.81282, -2.7305) (1.3016, 
 2.5593) (1.3774, -2.5240) (1.8197, 
 2.2632) (1.8865, -2.2142) (2.2642, 
 1.8744) (2.3192, -1.8138) (2.6169, 
 1.4089) (2.6580, -1.3390) (2.8635, 0.88570) (2.9937, 
 0.32621) (3.0022, -0.24664)
};
\draw[very thick] (0,0) circle (50mm);
\draw[-stealth, dashed,gray] (0,0) -- (32:50mm) 
node[below left, rotate=30, black]{$2+|\xi|$}; 
\foreach \x in {-4,...,-1}
\draw[shift={(\x,0)}] (0pt, 1pt) -- (0pt,-1pt) node[below] {\tiny $\x$};
\foreach \x in {1,...,4}
\draw[shift={(\x,0)}] (0pt, 1pt) -- (0pt,-1pt) node[below] {\tiny $\x$};
\end{tikzpicture}
\end{minipage}

\caption{\label{fig:1} Left: Zeros of the polar Jacobi polynomial
$P_{30}(z; 1/2, 2; 3 \exp(2\pi k I/30))$ for $k = 0, 1, \ldots, 29$.
Right: Zeros of the polar Jacobi polynomial
$P_{30}(z; \sqrt 3, \pi; 3 \exp(2\pi k I/23))$ for $k = 0, 1, \ldots, 22$.}
\end{figure}

In Figure \ref{fig:1}, we illustrate 
\rcrfi{on the} one hand how 
accurate Theorem \ref{thm:5} is, and 
\rcrfi{on} 
the other \rcrfi{hand}, we show the 
behavior of the zeros of the same polar Jacobi polynomial when the 
pole travels along a specific circle ({\rcrfi{{observe} 
 $I=\sqrt{-1}$}}).

\begin{figure}[H]
\begin{center}
\begin{minipage}{0.46\textwidth}
\begin{tikzpicture}[domain=-0.7:0.7,font=\sffamily, scale=0.85]
\draw[-stealth] (-3.3,0) -- (3.3,0) node[above] {$x$};
\draw[-stealth] (0,-3.3) -- (0,3.3) node[left] {$y$};
\draw plot [only marks, mark=*, mark options={fill=black},mark size=0.3pt]
coordinates{
(-0.0043832, -2.3480) (-2.8131, 
 0.59614) (0.13223, -2.4862) (-2.9279, 
 0.52644) (0.29249, -2.5828) (-3.0236, 
 0.42429) (0.46494, -2.6361) (-3.0915, 
 0.29650) (0.63881, -2.6464) (-3.1255, 
 0.15147) (0.80414, -2.6165) (-3.1217, -0.0013683) (0.95209, 
-2.5509) (-3.0786, -0.15202) (1.0752, -2.4560) (-2.9973, 
-0.29031) (1.1680, -2.3401) (-2.8810, -0.40628) (1.2271, 
-2.2123) (-2.7356, -0.49055) (1.2521, -2.0830) (-2.5696, 
-0.53484) (1.2458, -1.9623) (-2.3942, -0.53264) (1.2155, 
-1.8591) (-2.2240, -0.48047) (1.1720, -1.7784) (-2.0760, 
-0.38015) (1.1280, -1.7183) (-1.9674, -0.24139) (1.0922, 
-1.6697) (-1.9097, -0.082154) (-1.9045, 
 0.075975) (1.0651, -1.6199) (-1.9431, 
 0.21478) (1.0392, -1.5599) (-2.0116, 
 0.32327) (1.0031, -1.4873) (-2.0950, 
 0.39752) (0.94665, -1.4062) (0.86212, -1.3250) (-2.1796, 
 0.43924) (0.74561, -1.2553) (-2.2541, 
 0.45457) (0.59740, -1.2107) (-2.3104, 0.45336) (-2.3453, 
 0.44872) (0.42322, -1.2060) (-2.3630, 
 0.45531) (0.23642, -1.2561) (-2.3777, 
 0.48423) (0.060139, -1.3700) (-2.4103, 
 0.53408) (-0.076365, -1.5428) (-0.15089, -1.7520) (-2.4756, 
 0.58796) (-0.15811, -1.9695) (-2.5730, 
 0.62441) (-0.10532, -2.1729) (-2.6903, 0.62900)
 };
 \draw plot [only marks, mark=*, mark options={fill=white},mark size=0.8pt]
coordinates{
(0.2687, -1.104) (-1.480, -0.06404)
 };
\draw[very thick] (0,0) circle (30mm);
\draw[-stealth, dashed,gray] (0,0) -- (32:30mm) 
node[below left, rotate=30, black]{$2+|\xi|$}; 
\foreach \x in {-2,...,-1}
\draw[shift={(\x,0)}] (0pt, 1pt) -- (0pt,-1pt) node[below] {\tiny $\x$};
\foreach \x in {1,...,2}
\draw[shift={(\x,0)}] (0pt, 1pt) -- (0pt,-1pt) node[below] {\tiny $\x$};
\end{tikzpicture}
\end{minipage}
\hspace{3.5mm}
\begin{minipage}{0.5\textwidth}
\begin{tikzpicture}[domain=-0.7:0.7, font=\sffamily, scale=0.85]
\draw[-stealth] (-3.3,0) -- (3.3,0) node[above] {$x$};
\draw[-stealth] (0,-3.3) -- (0,3.3) node[left] {$y$};
\draw[gray] plot [only marks, mark=*, mark options={fill=black},mark size=0.3pt]
coordinates{(-1.5330, 
 0.60781) (-1.3978, -1.0052) (0.82640, -0.94776) (-1.6341, 
 0.67524) (-1.4419, -1.1513) (0.99351, -1.0770) (-1.7640, 
 0.69478) (-1.4284, -1.3130) (1.1745, -1.1336) (-1.8970, 
 0.66587) (-1.3642, -1.4677) (1.3479, -1.1311) (-2.0165, 
 0.59472) (-1.2579, -1.6018) (1.5010, -1.0812) (-2.1111, 
 0.48970) (-1.1184, -1.7058) (1.6251, -0.99502) (-2.1728, 
 0.36024) (-0.95484, -1.7731) (1.7143, -0.88329) (-2.1965, 
 0.21657) (-0.77720, -1.7990) (1.7649, -0.75724) (-2.1794, 
 0.069503) (-0.59591, -1.7808) (1.7755, -0.62839) (-2.1212, 
-0.069726) (-0.42197, -1.7179) (1.7478, -0.50857) (-2.0237, 
-0.18953) (-0.26708, -1.6119) (1.6864, -0.40975) (-1.8916, 
-0.27781) (-0.14414, -1.4670) (1.6005, -0.34344) (-1.7330, 
-0.32152) (-0.068242, -1.2927) (1.5059, -0.31866) (-1.5627, 
-0.30605) (-0.054636, -1.1112) (1.4265, -0.33465) (-1.4103, 
-0.21904) (-0.098430, -0.96478) (1.3826, -0.36924) (-1.3247, 
-0.074778) (-0.14451, -0.87775) (1.3649, -0.39262) (-1.3280, 
 0.066157) (-0.14638, -0.80085) (1.3482, -0.40254) (-1.3865, 
 0.15603) (-0.12985, -0.67715) (1.3256, -0.41728) (-1.4597, 
 0.18638) (-0.14995, -0.49646) (1.3144, -0.44728) (-1.5200, 
 0.16355) (-0.24830, -0.29209) (1.3331, -0.47345) (-1.5470, 
 0.094916) (-0.42613, -0.10561) (1.3688, -0.46843) (-1.5211, 
-0.020822) (-0.66861, 
 0.055774) (1.3944, -0.42904) (-1.4497, -0.22305) (-0.94252, 
 0.23969) (1.3924, -0.36726) (-1.4312, -0.46530) (-1.1318, 
 0.41322) (1.3542, -0.29854) (-1.4147, -0.65663) (-1.2753, 
 0.50166) (1.2767, -0.23912) (-1.3935, 
 0.52868) (-1.3721, -0.80132) (1.1612, -0.20647) (-1.4788, 
 0.51513) (-1.3093, -0.89611) (1.0146, -0.22102) (-1.5200, 
 0.48364) (-1.2481, -0.93273) (0.85469, -0.30827) (-1.5132, 
 0.47000) (-1.2336, -0.91689) (0.72887, -0.49152) (-1.4942, 
 0.51792) (-1.3049, -0.91681) (0.71663, -0.73834)
};
\draw[gray] plot [only marks, mark=*, mark options={fill=gray},mark size=0.6pt]
coordinates{
(-1.4375, -0.22728) (-0.95029, 
 0.41140) (-0.44193, -0.85664) (0.94818, -0.49429) (-1.5829, 
-0.28552) (-0.98084, 
 0.59491) (-0.42300, -1.0786) (1.1270, -0.60549) (-1.6816, 
-0.40985) (-1.1102, 
 0.71817) (-0.30811, -1.2634) (1.3049, -0.61851) (-1.7296, 
-0.56048) (-1.2681, 
 0.77219) (-0.14913, -1.3974) (1.4563, -0.56890) (-1.7291, 
-0.71915) (-1.4268, 
 0.76806) (0.032865, -1.4815) (1.5724, -0.47734) (-1.6831, 
-0.87279) (-1.5711, 
 0.71578) (0.22424, -1.5169) (1.6484, -0.35896) (-1.6903, 
 0.62480) (-1.5956, -1.0104) (0.41366, -1.5053) (1.6817, 
-0.22695) (-1.7768, 
 0.50478) (-1.4722, -1.1223) (0.59095, -1.4500) (1.6719, 
-0.093755) (-1.8249, 
 0.36594) (-1.3196, -1.2002) (0.74687, -1.3556) (1.6206, 
 0.028479) (-1.8311, 
 0.21918) (-1.1463, -1.2367) (0.87322, -1.2280) (1.5317, 
 0.12767) (-1.7941, 
 0.076043) (-0.96201, -1.2259) (0.96281, -1.0743) (1.4117, 
 0.19136) (-1.7145, -0.051077) (-0.77823, -1.1622) (1.0091, 
-0.90233) (1.2711, 
 0.20567) (-1.5953, -0.14826) (-0.61014, -1.0399) (1.0041, 
-0.71908) (1.1287, 
 0.15262) (-1.4430, -0.19782) (-0.48456, -0.85289) (0.92259, 
-0.53046) (1.0370, 
 0.0076205) (-1.2757, -0.17219) (-0.46152, -0.62501) (0.71810, 
-0.46596) (1.1158, -0.11157) (-1.1641, -0.047530) (-0.52241, 
-0.48904) (0.64723, -0.51812) (1.1577, -0.11212) (-1.1854, 
 0.074232) (-0.49848, -0.38537) (0.60818, -0.52749) (1.1723, 
-0.12028) (-1.2626, 
 0.10521) (-0.47852, -0.16838) (0.57720, -0.58153) (1.1959, 
-0.11538) (-1.3400, 0.055691) (-0.57625, 
 0.084865) (0.64259, -0.64090) (1.2011, -0.078684) (-1.4058, 
-0.052596) (-0.73470, 
 0.28558) (0.76919, -0.62936) (1.1588, -0.027290) (-1.4549, 
-0.18987) (-0.89941, 0.41810) (0.92638, -0.54998) (1.0463, 
 0.020966) (-1.4734, -0.33308) (-1.0558, 0.48273) (0.84724, 
 0.091182) (1.1095, -0.45658) (-1.4516, -0.46952) (-1.1966, 
 0.48795) (0.62489, 
 0.19033) (1.2463, -0.38105) (-1.3857, -0.58831) (-1.3130, 
 0.44318) (0.39132, 0.26648) (1.3212, -0.29364) (-1.3959, 
 0.35696) (-1.2738, -0.67784) (0.13627, 
 0.30396) (1.3428, -0.19884) (-1.4368, 
 0.23610) (-1.1115, -0.72632) (-0.14785, 
 0.29979) (1.3147, -0.11036) (-1.4317, 
 0.081755) (-0.88073, -0.73375) (-0.47808, 
 0.27322) (1.2398, -0.044887) (-1.4061, -0.10515) (-0.81188, 
 0.31608) (-0.59627, -0.76608) (1.1236, -0.023877) (-1.3887, 
-0.24879) (-0.99332, 
 0.36570) (-0.39390, -0.79665) (0.98083, -0.080341) (-1.3573, 
-0.28063) (-1.0389, 
 0.34652) (-0.33788, -0.76237) (0.87435, -0.26243)
 };
\draw[gray] plot [only marks, mark=*, mark options={fill=gray},mark size=1.2pt]
coordinates{
(-1.2539, 0.10993) (-1.0318, -0.43207) (-0.53025, 
 0.16402) (0.061173, -0.61042) (0.98023, -0.30137) (-1.4083, 
 0.14593) (-1.1151, -0.63887) (-0.52356, 
 0.48490) (0.13095, -0.88052) (1.1633, -0.38927) (-1.5556, 
 0.091899) (-1.0993, -0.84197) (-0.67048, 
 0.68358) (0.31012, -1.0460) (1.3272, -0.36410) (-1.6687, 
-0.011398) (-1.0265, -1.0230) (-0.84837, 
 0.79182) (0.50865, -1.1341) (1.4514, -0.28096) (-1.7414, 
-0.14458) (-1.0310, 
 0.83284) (-0.91027, -1.1757) (0.70694, -1.1615) (1.5321, 
-0.16403) (-1.7709, -0.29358) (-1.2043, 
 0.81888) (-0.76034, -1.2943) (0.89352, -1.1375) (1.5675, 
-0.029428) (-1.7564, -0.44608) (-1.3581, 
 0.75918) (-0.58590, -1.3735) (1.0592, -1.0699) (1.5577, 
 0.10938) (-1.6992, -0.59067) (-1.4838, 
 0.66255) (-0.39640, -1.4099) (1.1963, -0.96652) (1.5040, 
 0.24013) (-1.6022, -0.71661) (-1.5746, 
 0.53831) (-0.20157, -1.4012) (1.2987, -0.83620) (1.4097, 
 0.35130) (-1.6254, 
 0.39655) (-1.4705, -0.81380) (-0.011245, -1.3470) (1.2796, 
 0.43182) (1.3620, -0.68853) (-1.6327, 
 0.24823) (-1.3108, -0.87274) (0.16489, -1.2481) (1.1203, 
 0.47067) (1.3838, -0.53398) (-1.5952, 
 0.10531) (-1.1314, -0.88401) (0.31748, -1.1063) (0.93940, 
 0.45576) (1.3644, -0.38379) (-1.5130, -0.018701) (-0.94332, 
-0.83695) (0.43742, -0.92249) (0.74607, 
 0.37060) (1.3073, -0.25014) (-1.3886, -0.10679) (-0.76347, 
-0.71556) (0.51580, -0.68923) (0.55246, 
 0.18199) (1.2228, -0.14706) (-1.2294, -0.13030) (-0.64504, 
-0.49309) (0.31075, -0.22804) (0.63208, -0.33440) (1.1352, 
-0.091982) (-1.0986, -0.031876) (-0.69552, -0.30387) (0.14903, 
-0.38119) (0.77881, -0.26854) (1.0918, -0.084428) (-1.1387, 
 0.065616) (-0.65077, -0.16217) (0.093529, -0.41392) (0.82229, 
-0.28019) (1.0772, -0.071335) (-1.2304, 0.041789) (-0.69584, 
 0.11089) (0.15371, -0.52709) (0.88630, -0.16064) (1.0253, 
-0.12813) (-1.3195, -0.050205) (-0.83911, 
 0.30363) (0.37633, -0.59275) (0.66227, 
 0.022000) (1.1545, -0.16480) (-1.3816, -0.17127) (-0.99563, 
 0.39942) (0.41812, 
 0.20644) (0.62270, -0.66236) (1.2310, -0.098985) (-1.4022, 
-0.30599) (-1.1462, 0.42113) (0.19895, 
 0.36943) (0.81743, -0.68423) (1.2575, -0.0042215) (-1.3769, 
-0.44404) (-1.2797, 0.38583) (-0.022437, 
 0.48818) (0.97952, -0.64794) (1.2341, 0.099112) (-1.3873, 
 0.30668) (-1.3060, -0.57566) (-0.24879, 
 0.56373) (1.1116, -0.56673) (1.1605, 0.19660) (-1.4619, 
 0.19583) (-1.1930, -0.69138) (-0.47324, 0.59783) (1.0376, 
 0.27704) (1.2115, -0.45471) (-1.4982, 
 0.065902) (-1.0441, -0.78136) (-0.68517, 0.59210) (0.86851, 
 0.33184) (1.2754, -0.32734) (-1.4913, -0.069506) (-0.87286, 
 0.54757) (-0.86753, -0.83396) (0.65905, 
 0.35249) (1.2981, -0.20047) (-1.4366, -0.19581) (-1.0259, 
 0.46472) (-0.67018, -0.83383) (0.41400, 
 0.32699) (1.2750, -0.088829) (-1.3271, -0.29795) (-1.1371, 
 0.34428) (-0.44945, -0.75788) (0.12537, 
 0.23828) (1.2047, -0.0088482) (-1.2091, 
 0.19134) (-1.1491, -0.36857) (-0.28479, 
 0.10336) (-0.13578, -0.60321) (1.0907, 0.013918) (-1.2455, 
 0.062869) (-0.95445, -0.42644) (-0.64244, 
 0.14077) (0.13220, -0.56385) (0.95748, -0.075336)
};
 \draw plot [only marks, mark=+, mark options={fill=gray},mark size=1.3pt]
coordinates{
 (-1.0864, -0.030899) (-0.71543, -0.27618) (0.16859, -0.39328) (0.92800, -0.23309)
}; 
\draw plot [only marks, mark=x, mark options={fill=gray},mark size=1.3pt]
coordinates{
(-1.0518, -0.020481) (-0.82070, -0.16939) (-0.21802, -0.27836) (0.48301, -0.28089) (0.96208, -0.14247)
}; 
\draw plot [only marks, mark=*, mark options={fill=white},mark size=0.7pt]
coordinates{
(-1.1715, -0.049823) (-0.48065, -0.51329) (0.82386, -0.44574)
};

\draw[very thick] (0,0) circle (30mm);
\draw[-stealth, dashed,gray] (0,0) -- (32:30mm) 
node[below left, rotate=30, black]{$2+|\xi|$}; 
\foreach \x in {-2,...,-1}
\draw[shift={(\x,0)}] (0pt, 1pt) -- (0pt,-1pt) node[below] {\tiny $\x$};
\foreach \x in {1,...,2}
\draw[shift={(\x,0)}] (0pt, 1pt) -- (0pt,-1pt) node[below] {\tiny $\x$};
\end{tikzpicture}
\end{minipage}
\end{center}
\caption{\label{fig:3} Left: Zeros of the polar Jacobi polynomial
$P_{2}(z; -1/2+I, -1.45-I/2; \exp(2\pi k I/30))$ for $k = 0, 1, \ldots, 29$ (dots) and zeros of the Jacobi polynomials $P^{( -1/2+I, -1.45-I/2)}_{2}(z)$ (circles).
Right: Zeros of the polar Jacobi polynomial
$P_{n}(z; -1/2+I, -1.45-I/2; \exp(2\pi k I/30))$ for $k = 0, 1, \ldots, 29$ (gray dots) and zeros of the Jacobi polynomials $P^{( -1/2+I, -1.45-I/2)}_{n}(z)$ ($+$,$\times$, and circles) for $k=0, 1, \ldots, 29$, $n=3,4,5$.}
\end{figure}
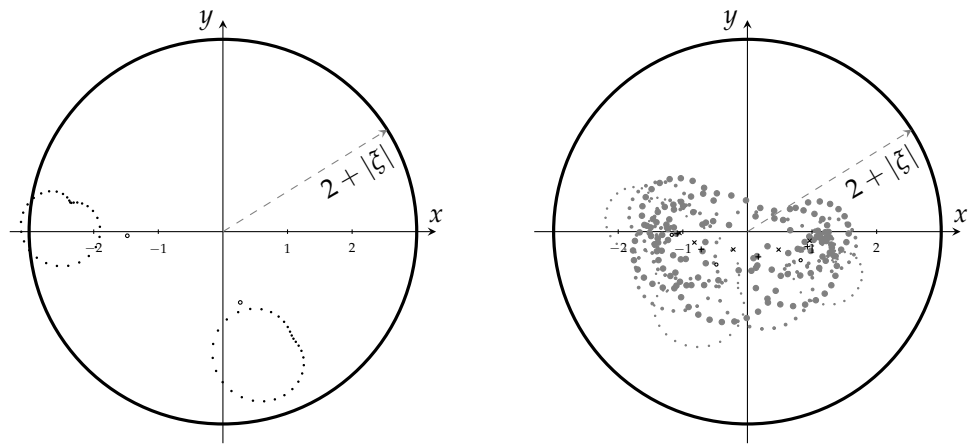

In Figure \ref{fig:3}, we illustrate an example of 
Jacobi polar polynomials where the parameters $\Re \alpha \le -1$ 
or $\Re \beta\le -1$; therefore, the zeros of the Jacobi polynomial 
can move away from the interval $[-1, 1]$ in a somewhat uncontrolled 
way. Therefore, Theorem \ref{thm:5} cannot be applied in such a case. 
However, observe that in the considered example $-2 <\Re(\alpha +\beta)
 = -1.95 < -1$.
 
\rcrfi{The next theorem gives the location of the zeros of the polar Jacobi 
polynomial of degree $n$ and its multiplicity, or equivalently the 
location of source points and their corresponding strength.}

\begin{Theorem}\label{thm:3.4}
For any $\Re \alpha,\Re \beta>-1$ and $\xi\in \mathbb C$.
The following statements hold:
\begin{enumerate} 
\item If $\zeta\in \mathbb C^*$ is a zero of 
$P_n(z;\alpha,\beta;\xi)$, then $z=-\zeta$ is a zero of
$P_n(z;\beta,\alpha;-\xi)$.
\item If $\zeta\in \mathbb C^*$ is a zero of 
$P^{(\alpha,\beta)}_n(z)$, then $\zeta$ is a zero of 
\rcrfi{$P_n(z)$}.
\item The zeros of 
\rcrfi{$P_n(z)$} have multiplicity of at most 2 and their 
multiple zeros are located on $[-1,1]$.
\item All the zeros of 
\rcrfi{$P_n(z)$} are located on the curve
\begin{equation} \label{eq:34}
{\mathscr Z}_n(\xi)=\left\{z\in \mathbb C\ :\ P_{n+1}^{(\alpha-1,\beta-1)}(z)=
P_{n+1}^{(\alpha-1,\beta-1)}(\xi)\right\}\setminus\{\xi\}.
\end{equation}
\end{enumerate}
\end{Theorem}
\begin{proof} 
The first statement holds true due to {equation} \eqref{eq:34}, 
the second statement holds true due to {equation} \eqref{eq:relpoljac},
and the fourth statement holds true due to {equation} \eqref{eq:RJPoJ2}.

\rcrfi{Suppose that} 
$\omega$ is a zero $P_n$ of multiplicity greater 
than two; then, by  \eqref{eq:relpoljac}, $\omega$ is a zero 
of $P_n^{(\alpha,\beta)}$ and also a zero of 
$\left(P_n^{(\alpha,\beta)}\right)'$. Thus, $\omega$ is a 
zero of multiplicity 2 of $P_n^{(\alpha,\beta)}$. This is a 
contradiction since the zeros of the Jacobi polynomials 
are all simple. Therefore, statement 3 holds true.
\end{proof} 
\begin{Remark}~
\begin{itemize} 
\item Observe that the zeros of $P_n$ do not have to be simple. 
Let $\xi_+=(1+2\sqrt 6)/5$ or $\xi_-=(1-2\sqrt 6)/5$; then, the 
polar polynomial of degree two
$
P_2(z;0,1,\xi_+)=\left(z-\frac {1-\sqrt 6}5\right)^2$, 
or $P_2(z;0,1,\xi_-)=\left(z-\frac {1+\sqrt 6}5\right)^2$.
\item When the parameters are 
\rcrfi{not} 
standard, i.e.,  $\Re \alpha<-1$ or $\Re \beta<-1$ then, 
by Corollary \ref{cor:2}, statement 3 of 
Theorem \ref{thm:3.4} is no longer true. 
For example, if $\alpha=-4$, $\beta=1>0$, 
and $n=5$, then $P_5(z;-4,1,1)=(z-1)^4 (z-5/7)$.
\end{itemize}
\end{Remark}

We can establish the following result concerning the 
boundedness of the zeros of the polar polynomials.
\begin{Lemma} \label{lem:3}
Given $\xi\in \mathbb C$, let us define the two numbers 
$\Delta_\xi:=\sup\{|\xi-z| : z\in [-1,1]\}$ and
$\delta_\xi:=\inf\{|\xi-z| : z\in [-1,1]\}$. Then, the following can be stated: 
\begin{enumerate} 
\item All zeros of the polar Jacobi polynomials with 
pole $\xi$ are contained in $|z|\le \Delta_\xi+1$.
\item If $\delta_\xi>1$, the zeros of the polar Jacobi 
polynomials with pole $\xi$ are simple and contained 
in the exterior of the ellipse $|z+1|+|z-1|=2\alpha$, 
where $1<\alpha<\delta_\xi$.
\end{enumerate}
\end{Lemma}
\begin{proof} 
By {equation} \eqref{eq:RJPoJ2}, the zeros of $P_n(z)$ are 
located in ${\mathscr Z}_n(\xi)$. Since 
$\left|P_{n+1}^{(\alpha-1,\beta-1)}(\xi)\right|<\Delta_\xi^{n+1}$,
they are contained in the interior of the set 
$\left|P_{n+1}^{(\alpha-1,\beta-1)}(z)\right|=\Delta_\xi^{n+1}$.
It is known 
\rcrfi{that} the zeros of $P_{n+1}^{(\alpha-1,\beta-1)}(z)$,
namely $x_{n+1,k}$, satisfy $|x_{n+1,k}|\le 1$. 
Therefore, for any $t\in \mathbb C$, such that 
$|t|>1+\Delta_\xi$, we have 
\[
\left|P_{n+1}^{(\alpha-1,\beta-1)}(z)\right|=\prod_{k=0}^n 
|z-x_{n+1,k}|\ge \prod_{k=0}^n \big||z|-|x_{n+1,k}|\big|>
\Delta_\xi^{n+1}.
\]
Hence, the first statement holds. 

\rcrfi{Concerning} 
the second statement, let $z$ be such that 
$|z+1|+|z-1|=2\alpha$. From the well-known 
arithmetic--geometric mean inequality, we have
\[
\left|P_{n+1}^{(\alpha-1,\beta-1)}(z)\right|\le 
\left(\frac 1{n+1}\sum_{k=0}^n |z-x_{n+1,k}|\right)^{n+1}
<\alpha^{n+1}.
\]
{If} $\omega$ is a zero of $P_n$, 
\rcrfi{we obtain, in view of} 
{equation} \eqref{eq:34},
\[
\left|P_{n+1}^{(\alpha-1,\beta-1)}(\omega)\right|=
\left|P_{n+1}^{(\alpha-1,\beta-1)}(\xi)\right|=
\prod_{k=0}^n |\xi-x_{n+1,k}|>\delta_\xi^{n+1}>\alpha^{n+1}.
\]
{Therefore}, the result holds.
\end{proof}
The last result is about the asymptotic behavior of the zeros of the 
polar Jacobi polynomials.
\begin{Theorem}[Theorem 22 in \cite{MR3032627}]\label{thm:7}
The accumulation points of zeros of $(P_n)$ are located on the set 
${\mathscr Z}(\xi)\cup [-1,1]$, where ${\mathscr Z}(\xi)$ is the
ellipse
\vspace{-12pt}{}
\begin{adjustwidth}{-\extralength}{0cm}
\centering 
\[
{\mathscr Z}(\xi)=\{z\in \mathbb C\, : \, 
z=\cosh(\log |\varphi(\xi)|+i\theta),\ 0\le \theta<2\pi\}
=\left\{z\in \mathbb C\, : \, \left|z+\sqrt{z^2-1}\right|=
\left|\varphi(\xi)\right|\right\},
\]
\end{adjustwidth}
where $\varphi(z)=z+\sqrt{z^2-1}$.
\end{Theorem}


\dataavailability{The original contributions presented in the study are included in the article, further inquiries can be directed to the corresponding author.}

\conflictsofinterest{The authors declare no conflicts of interest.} 

\begin{adjustwidth}{-\extralength}{0cm}

\reftitle{References}

\PublishersNote{}
\end{adjustwidth}
\end{document}